\documentclass[english]{article}
\usepackage[T1]{fontenc}
\usepackage[latin9]{inputenc}
\usepackage{geometry}
\usepackage{babel}
\usepackage{amsthm}
\usepackage{amsmath}
\usepackage{amssymb}
\usepackage[unicode=true,pdfusetitle,
 bookmarks=true,bookmarksnumbered=false,bookmarksopen=false,
 breaklinks=false,pdfborder={0 0 1},backref=section,colorlinks=false]
 {hyperref}

\makeatletter
\theoremstyle{plain}
\newtheorem{thm}{\protect\theoremname}[section]
  \theoremstyle{definition}
  \newtheorem{defn}[thm]{\protect\definitionname}
  \theoremstyle{plain}
  \newtheorem{fact}[thm]{\protect\factname}
  \theoremstyle{definition}
  \newtheorem{problem}[thm]{\protect\problemname}
  \theoremstyle{plain}
  \newtheorem{lem}[thm]{\protect\lemmaname}
  \theoremstyle{remark}
  \newtheorem{rem}[thm]{\protect\remarkname}
  \theoremstyle{plain}
  \newtheorem{cor}[thm]{\protect\corollaryname}
  \theoremstyle{remark}
  \newtheorem*{claim*}{\protect\claimname}

\makeatother

  \providecommand{\claimname}{Claim}
  \providecommand{\corollaryname}{Corollary}
  \providecommand{\definitionname}{Definition}
  \providecommand{\factname}{Fact}
  \providecommand{\lemmaname}{Lemma}
  \providecommand{\problemname}{Problem}
  \providecommand{\remarkname}{Remark}
\providecommand{\theoremname}{Theorem}

\begin{document}

\title{\noindent \textbf{\Large{}Indiscernible arrays and rational functions
with algebraic constraint}}

\author{\noindent ELAD LEVI}

\maketitle
\noindent \begin{flushleft}
{\footnotesize{}~~~Abstract. Let $k$ be an algebraically closed
field of characteristic zero and $P(x,y)\in k[x,y]$ be a polynomial
which dep-}\\
{\footnotesize{}~~~ends on all its variables. P has an algebraic
constraint if the set $\mbox{\ensuremath{\{(P(a,b),(P(a',b'),P(a',b),P(a,b')\,|\,a,a',b,b'\in k\}}}$
}\\
{\footnotesize{}~~~does not have the maximal Zariski-dimension.
Tao proved that if $P$ has an algebraic constraint then it can be
decom-}\\
{\footnotesize{}~~~posed: there exists $Q,F,G\in k[x]$ such that
$P(x_{1},x_{2})=Q(F(x_{1})+G(x_{2})),$ or $P(x_{1},x_{2})=Q(F(x_{1})\cdot G(x_{2}))$.
In this }\\
{\footnotesize{}~~~paper we give an answer to a question raised
by Hrushovski and Zilber regarding 3-dimensional indiscernible arrays
in}\\
{\footnotesize{}~~~stable theories. As an application of this result
we find a decomposition of rational functions in three variables which}\\
{\footnotesize{}~~~has an algebraic constraint.}
\par\end{flushleft}{\footnotesize \par}

\section{Introduction}

\subsubsection*{Rational functions with an algebraic constraint}

Let $k$ be an algebraically closed field of characteristic zero,
and let $r(x_{1},...,x_{n})\in k(x_{1},..,x_{n})$ be a rational function.
For every $(v_{1}^{0},..,v_{n}^{0},v_{1}^{1},..,v_{n}^{1})\in k^{2n}$
let 
\begin{eqnarray*}
 &  & L_{r}(v_{1}^{0},..,v_{n}^{0},v_{1}^{1},..,v_{n}^{1})=(r(v_{1}^{0},..,v_{n}^{0}),r(v_{1}^{1},v_{2}^{0},..,v_{n}^{0}),..,r(v_{1}^{1},.,v_{n}^{1}))\in k^{2^{n}}.
\end{eqnarray*}
$r(x_{1},..,x_{n})$ has an algebraic constraint if the set $\mbox{\ensuremath{A_{r}=\{L_{r}(v_{1}^{0},..,v_{n}^{0},v_{1}^{1},..,v_{n}^{1})\,|\,\ v_{1}^{0},..,v_{n}^{1}\in k\}}}$
has Zariski dimension less then $2n$. We will say that $r(x_{1},..,x_{n})$
is non-degenerate, if $r(v_{1},..,v_{n})\notin k(v_{i_{1}},..,v_{i_{n-1}})$
for every algebraic independent set $\{v_{1},..,v_{n}\}\subset k$
and every $i_{1},..,i_{n-1}\in\{1,..,n\}$.\\
\\
Tao proved in {[}1{]}(Theorem 41) that in case $r(x_{1},x_{2})\in k[x_{1},x_{2}]$
is a non-degenerate polynomial with an algebraic constraint, then
there exists $Q,F,G\in k[x]$ such that either $r(x_{1},x_{2})=Q(F(x_{1})+G(x_{2})),$
or $r(x_{1},x_{2})=Q(F(x_{1})\cdot G(x_{2}))$. Hrushovski {[}2{]}
observed that this theorem can also be proved by the \\
celebrated Group Configuration theorem (see {[}3{]} for more information
about this theorem). He \\
managed to prove that if $\mbox{\ensuremath{r(x_{1},x_{2})\in k(x_{1},x_{2})}}$
is a non-degenerate rational function with an algebraic constraint,
then there exist a one dimensional algebraic group $(G,*)$, rational
maps $\mbox{\ensuremath{\overline{r}\in k(x_{1},x_{2}),\,\,\pi\in k(x_{1})}}$,
and $f,g,h,q:\,k\rightarrow G$ such that: $\pi\circ\overline{r}=r$
and $q\circ\overline{r}(u_{1},u_{2})=g(u_{1})*h(u_{2})$.\\
\\
In this paper we focus on rational functions in three variables that
are non-degenerate and have algebraic constraint. In section 6 we
prove the following result:
\begin{thm}
Let $P(x,y,z)$ be a rational function over an algebraically closed
field $k$ of characteristic zero. If $P$ is non-degenerate and has
an algebraic constraint, then one of the following holds:
\begin{enumerate}
\item There exists $n\in\mathbb{N}$, rational functions $r_{1},r_{2},r_{3},q,\pi$
and $\overline{P}\in k(x,y,z)$ such that for some $i\neq j\neq l\in\{1,2,3\}$:
$q\circ\overline{P}(x_{1},x_{2},x_{3})=r_{i}(x_{i})(r_{j}(x_{j})+r(x_{l}))^{n}\mbox{ and }\pi\circ\overline{P}(x,y,z)=P(x,y,z).$ 
\item There exists rational functions $r_{1},r_{2},r_{3},q,\pi\in k(x)$
and $\overline{P}\in k(x,y,z)$ such that:\\
$q\circ\overline{P}(x,y,z)=\frac{r_{1}(x)+r_{2}(y)}{r_{2}(y)+r_{3}(z)}\mbox{ and }\pi\circ\overline{P}(x,y,z)=P(x,y,z).$
\item There exists a one dimension algebraic group $(G,*)$, rational maps
\textup{$r_{1},r_{2},r_{3},q,\pi$ and $\overline{P}$ such that:
$q\circ\overline{P}(x,y,z)=r_{1}(x)*r_{2}(y)*r_{3}(z)$ and }$\pi\circ\overline{P}(x,y,z)=P(x,y,z)$. 
\end{enumerate}
\end{thm}
In particular, we use this theorem to find a family of rational functions
(definition \ref{Def:2-dec}) that has an algebraic constraint. Our
proof of this theorem is based on Hrushovski's ideas; these lead us
to the main topic of this paper.

\subsubsection*{Indiscernible arrays}

Let $M$ be a model, an array of elements: $F=(f_{i,j})_{i,j\in\mathbb{N}}\subset M$
is an indiscernible array if for every $n\in\mathbb{N}$ and every
natural numbers $\mbox{\ensuremath{i_{1}<i_{2}<..<i_{n},\,\,\ j_{1}<j_{2}<..<j_{n}}}$
we have: 
\begin{eqnarray*}
 &  & tp(f_{t,k}\,|\,1\leq t,k\leq n)=tp(f_{i_{t},j_{k}}\,|\,1\leq t,k\leq n).
\end{eqnarray*}
In the case where M is strongly minimal we can use Morely rank to
define the rank function of the indiscernible array: $\mbox{\ensuremath{\alpha_{F}(m,n)=RM(f_{i,j}\,|\,1\leq i\leq m,\,\,1\leq j\leq n)}}$.\\
\\
Hrushovski and Zilber used {[}4{]}(6.3) the Group Configuration theorem
to prove that if $\mbox{\ensuremath{\alpha_{F}(m,n)=\lambda m+n-\lambda}}$
then $F$ is induced by a definable action of definable group of Morely
rank $\lambda$. Hrushovski and Zilber made a prediction {[}4{]}(6.5)
for the possibilities of a 3-dimension indiscernible array with rank
function $\mbox{\ensuremath{\alpha_{F}(m,n.p)=m+n+p-2}}$. In this
paper we show another possibility that was not considered by Hrushovski
and Zilber (section 3). However, we prove (sections 4-5) that this
possibility is the only exception. In the next section we give the
basic definitions, and describe precisely the main problem that we
wish to answer.

\subsubsection*{Acknowledgments}

The author is greatly indebted to Ehud Hrushovski for his guidance
and for sharing his ideas, and Yatir Halevi for a careful reading
of an earlier draft. The research leading to this paper has been partially
supported by the European Research Council under the European Union's
Seventh Framework Programme (FP7/2007-2013)/ERC Grant Agreement No.
291111

\section{Preliminaries}

Let T be a stable theory, we work inside a monster model $\mathcal{U}\vDash T$
(a big saturated model) in which $\mbox{\ensuremath{acl(\emptyset)=dcl(\emptyset)}}$.
\\
We denote the Morley rank by $RM$, and use $\downarrow$ to denote
the non-forking independence relation.  
\begin{defn}
\end{defn}
\begin{enumerate}
\item $F=(f_{i^{1},...,i^{n}})_{i^{1},..,i^{n}\in\mathbb{N}}\subset\mathcal{U}$
is an n indiscernible array if $RM(f_{i^{1},...,i^{n}})=1$ and for
every $m\in\mathbb{N}$, every $j_{1}^{t}<j_{2}^{t}<...<j_{m}^{t}$
(where $1\leq t\leq n$); $\mbox{tp(\ensuremath{f_{j_{i_{1}}^{1},..,j_{i_{n}}^{n}}}\,|\,1\ensuremath{\leq i_{1}},..,\ensuremath{i_{n}\leq}m)=tp(\ensuremath{f_{i_{1},..,i_{n}}}\,|\,1\ensuremath{\leq i_{1}},..,\ensuremath{i_{n}\leq}m)}$.
\item For every n indiscernible array $F=(f_{i^{1},...,i^{n}})_{i^{1},..,i^{n}\in\mathbb{N}}$
we define the rank function:
\begin{eqnarray*}
 &  & \alpha_{F}(m_{1},..,m_{n})=RM(f_{i^{1},..,i^{n}}\,|\,1\leq i^{j}\leq m_{j},\,\,1\leq j\leq n).
\end{eqnarray*}

\item We say that two n indiscernible arrays $F=(f_{i^{1},...,i^{n}})_{i^{1},..,i^{n}\in\mathbb{N}}$
and $F'=(f'_{i^{1},...,i^{n}})_{i^{1},..,i^{n}\in\mathbb{N}}$ are
algebraically equivalent, if for every $i^{1},..,i^{n}\in\mathbb{N}$
we have: $\mbox{\ensuremath{acl(f_{i^{1},..,i^{n}})=acl(f'_{i^{1},..,i^{n}})} }$.
\end{enumerate}
In this paper we only focus on $2$ and 3 indiscernible arrays. We
call 2 indiscernible arrays simply indiscernible arrays.\\
\\
As mentioned in {[}4{]}(6.2), there is a polynomial $p\in\mathbb{Z}[x,y,z]$
such that $\mbox{\ensuremath{p(m,n,r)=\alpha_{F}(m,n,r)}}$ for large
enough $m,n,r$. A special case is when $F=(f_{i,j})$ is an indiscernible
array with $\mbox{\ensuremath{\alpha_{F}(m,n)=\lambda m+n-\lambda}}$,
in this case Hrushovski and Zilber used the celebrated Group Configuration
theorem in order to prove the following useful theorem:
\begin{fact}
\label{fct:2.2} $[4](6.3)$ Let $F=(f_{i,j})$ be an indiscernible
array.
\begin{enumerate}
\item If $\alpha_{F}(m,n)=m+n-1$ then there exists a definable one dimensional
group $(G,\cdot)$ and elements $\{b_{i}\}_{i\in\mathbb{N}},\{a_{j}\}_{j\in\mathbb{N}}\subset G$
such that $acl(f_{i,j})=acl(a_{i}\cdot b_{j})$ for every $i,j\in\mathbb{N}$.
\item If $\alpha_{F}(m,n)=2m+n-2$ then there exists a definable field $(\mathbb{F},+,\cdot)$
and elements $\mbox{\{\ensuremath{b_{i}\}_{i\in\mathbb{N}}},\{\ensuremath{a_{j}\}_{j\in\mathbb{N}}},\{\ensuremath{c_{k}\}_{k\in\mathbb{N}}\subset\mathbb{F}}}$
such that $acl(f_{i,j})=acl(a_{i}\cdot c_{j}+b_{j})$ for every $i,j\in\mathbb{N}$.
\end{enumerate}
\end{fact}
\begin{defn}
\label{Definition}Let $F=(f_{i,j,k})$ be a 3 indiscernible array.
\begin{enumerate}
\item We say that F has the group form if there exists a definable one dimensional
group $(G,\cdot)$ and elements $\{b_{i}\}_{i\in\mathbb{N}},\{a_{j}\}_{j\in\mathbb{N}},\{c_{k}\}_{k\in\mathbb{N}}\subset G$
such that $acl(f_{i,j,k})=acl(a_{i}\cdot b_{j}\cdot c_{k})$.
\item We say that $F$ has the field form if there exists a definable field
$(\mathbb{F},+,\cdot)$ , elements $\mbox{\{\ensuremath{a_{i}^{1}\}_{i\in\mathbb{N}}},\{\ensuremath{a_{j}^{2}\}_{j\in\mathbb{N}}},\{\ensuremath{a_{k}^{3}\}_{k\in\mathbb{N}}\subset\mathbb{F}}}$,
and some permutation $\sigma\in S_{3}$ such that $acl(f_{i_{1},i_{2},i_{3}})=acl(a_{i_{\sigma(1)}}^{\sigma(1)}\cdot(a_{i_{\sigma(2)}}^{\sigma(2)}+a_{i_{\sigma(3)}}^{\sigma(3)}))$
for every $i_{1},i_{2},i_{3}\in\mathbb{N}$.
\item We  say that $F$ has the twisted form if there exists a definable
field $(\mathbb{F},+,\cdot)$ and an independent set $\{a_{i}^{1},a_{j}^{2},a_{k}^{3}\,|\,i,j,k\in\mathbb{N}\}\subset\mathbb{F}$
such that: $\mbox{acl(\ensuremath{f_{i_{1},i_{2},i_{3}}})=acl(\ensuremath{\frac{a_{i}^{1}+a_{j}^{2}}{a_{k}^{3}+a_{j}^{2}}}) }$
for every $i,j,k\in\mathbb{N}$.
\end{enumerate}
\end{defn}
Hrushovski and Zilber asked {[}4{]}(6.5) the following question:
\begin{problem}
\label{conjecture 1}Let $F=(f_{i,j,k})$ be a 3 indiscernible array
with rank function $\alpha_{F}(m,n,k)=m+n+k-2$. Are the following
the only possibilities?
\begin{enumerate}
\item F has the group form.
\item F has the field form.
\end{enumerate}
\end{problem}
In the next section we give a negative answer to this problem, in
other words, there is another case which was not considered- the twisted
form. However, the twisted form is the only exception. The main purpose
of this paper is to prove the following theorem:
\begin{thm}
\label{thm:main}Let $F=(f_{i,j,k})$ be a 3 indiscernible array and
suppose that $\alpha_{F}(m,n,k)=m+n+k-2$, then there are only three
possibilities:
\begin{enumerate}
\item F has the group form.
\item F has the field form.
\item F has the twisted form.
\end{enumerate}
\end{thm}

\section{The twisted form- A counterexample}

In this section we will see that that the twisted form is really a
counterexample. We start with some basic fact regarding indiscernible
arrays:
\begin{lem}
\textup{\label{lem:Counter1}} Indiscernible arrays have the following
properties:
\begin{enumerate}
\item Let $F=(f_{i,j,k})$ be a 3 indiscernible array array, then $\mbox{\ensuremath{\alpha_{F}}(m,n,k)=m+n+k-2}$
iff 
\begin{eqnarray*}
 &  & \alpha_{F}(1,m,n)=\alpha_{F}(m,1,n)=\alpha_{F}(m,n,1)=m+n-1.
\end{eqnarray*}

\item Let $F=(f_{i,j})$ be an indiscernible array, then $\alpha_{F}(m,n)=m+n-1$
iff $\mbox{\ensuremath{RM(f_{1,1},f_{1,2},f_{2,1},f_{2,2})=3}}$.\end{enumerate}
\begin{proof}
$ $
\begin{enumerate}
\item The left to right direction is trivial. For the other direction, assume
\begin{eqnarray*}
 &  & \alpha_{F}(1,m,n)=\alpha_{F}(m,1,n)=\alpha_{F}(m,n,1)=m+n-1.
\end{eqnarray*}
We have: $\alpha_{F}(2,1,2)=3$, thus by indiscernibility: $f_{i,1,k}\in acl(f_{1,1,1},f_{1,1,k},f_{i,1,1})$
for every $2\leq k$. We also know from $\alpha_{F}(1,m,n)=m+n-1$
and indiscernibility that: 
\begin{eqnarray*}
 &  & f_{i,j,k}\in acl(f_{i,j,1},f_{i,1,k},f_{i,1,1}).
\end{eqnarray*}
In conclusion, for every $m,n,p\in\mathbb{N}$ we have:
\begin{eqnarray*}
 &  & \alpha(m,n,p)=RM((f_{i,j,k}\,|\,1\leq i\leq m,\quad1\leq j\leq n,\quad1\leq k\leq p))=\\
 &  & RM(f_{i,j,1}\,|\,1\leq i\leq m,\quad1\leq j\leq n)+\\
 &  & +RM((f_{i,j,k}\,|\,1\leq i\leq m,\quad1\leq j\leq n,\quad2\leq k\leq p)/(f_{i,j,1}\,|\,1\leq i\leq m,\quad1\leq j\leq n))=\\
 &  & m+n-1+RM((f_{1,1,k}\,|\,2\leq k\leq p)/(f_{i,j,1}\,|\,1\leq i\leq m,\quad1\leq j\leq n))=\\
 &  & m+n+p-2
\end{eqnarray*}

\item Again, we only need to prove the direction from right to left: If
$\mbox{RM(\ensuremath{f_{1,1}},\ensuremath{f_{1,2}},\ensuremath{f_{2,1}},\ensuremath{f_{2,2}})=3}$
then by indiscernibility for every $i,j\in\mathbb{N}$: $RM(f_{1,i},f_{1,1},f_{j,1},f_{i,j})=3$,
thus $f_{i,j}\in acl(f_{1,i},f_{1,1},f_{j,1})$ and we get for every
$m,n\in\mathbb{N}$:
\begin{eqnarray*}
 &  & \alpha(m,n)=RM((f_{i,j}\,|\,1\leq i\leq m,\quad1\leq j\leq n))=\\
 &  & RM(f_{i,1}\,|\,1\leq i\leq m)+RM((f_{i,j}\,|\,1\leq i\leq m,\quad2\leq j\leq n)/(f_{i,1}\,|\,1\leq i\leq m))=\\
 &  & m+RM((f_{1,j}\,|\,2\leq j\leq n)/(f_{i,1}\,|\,1\leq i\leq m))=m+n-1
\end{eqnarray*}

\end{enumerate}
\end{proof}
\end{lem}
The following lemma proves that every 3 indiscernible array which
has the twisted form satisfies the assumptions of problem \ref{conjecture 1}:
\begin{lem}
\label{lem: equal 2}Let $F=(f_{i,j,k})$ be a 3 indiscernible array,
if $F$ has the twisted form then\\
$\alpha_{F}(m,n,k)=m+n+k-2$.
\begin{proof}
Let $(\mathbb{F},+,\cdot)$ be a field, $\{v_{i}\}_{i\in\mathbb{N}},\{u_{j}\}_{j\in\mathbb{N}},\{w_{k}\}_{k\in\mathbb{N}}\subset\mathbb{F}$
such that: $\mbox{acl(\ensuremath{f_{i,j,k}})=acl(\ensuremath{\frac{u_{i}+v_{j}}{v_{j}+w_{k}}})}$
for every $i,j,k\in\mathbb{N}$ and let $t_{i,j,k}=\frac{u_{i}+v_{j}}{v_{j}+w_{k}}$.
By Lemma \ref{lem:Counter1}(1) it will suffice to show that 
\begin{eqnarray*}
 &  & \alpha_{F}(1,m,n)=\alpha_{F}(m,1,n)=\alpha_{F}(m,n,1)=m+n-1.
\end{eqnarray*}
 Notice that:
\begin{enumerate}
\item $\frac{t_{2,1,1}}{t_{1,1,1}}=\frac{\frac{u_{2}+v_{1}}{v_{1}+w_{1}}}{\frac{u_{1}+v_{1}}{v_{1}+w_{1}}}=\frac{u_{2}+v_{1}}{u_{1}+v_{1}}=\frac{\frac{u_{2}+v_{1}}{v_{1}+w_{2}}}{\frac{u_{1}+v_{1}}{v_{1}+w_{2}}}=\frac{t_{2,1,2}}{t_{1,1,2}}$.
\item $\frac{t_{2,1,1}-1}{t_{1,1,1}-1}=\frac{\frac{u_{2}-w_{1}}{v_{1}+w_{1}}}{\frac{u_{1}-w_{1}}{v_{1}+w_{1}}}=\frac{u_{2}-w_{1}}{u_{1}-w_{1}}=\frac{\frac{u_{2}-w_{1}}{v_{2}+w_{1}}}{\frac{u_{1}-w_{1}}{v_{2}+w_{1}}}=\frac{t_{2,2,1}-1}{t_{1,2,1}-1}$.
\item $\frac{t_{1,1,2}^{-1}-1}{t_{1,1,1}^{-1}-1}=\frac{\frac{v_{1}+w_{2}}{u_{1+}v_{1}}-1}{\frac{v_{1}+w_{1}}{u_{1+}v_{1}}-1}=\frac{w_{2}-u_{1}}{w_{1}-u_{1}}=\frac{\frac{v_{2}+w_{2}}{u_{1+}v_{2}}-1}{\frac{v_{2}+w_{1}}{u_{1+}v_{2}}-1}=\frac{t_{1,2,2}^{-1}-1}{t_{1,2,1}^{-1}-1}$.
\end{enumerate}

Combining lemma \ref{lem:Counter1}(2) with (1)+(2)+(3) we get: $\mbox{\ensuremath{\alpha_{F}(1,m,n)=\alpha_{F}(m,1,n)=\alpha_{F}(m,n,1)=m+n-1}}$
as desired.

\end{proof}
\end{lem}
In order to show  that the twisted form is a counterexample to Problem
\ref{conjecture 1}, we must find a dividing line which will separate
the forms in Definition \ref{Definition}.
\begin{defn}
\label{def:frames}Let $F=(f_{i,j,k})$ be a 3 indiscernible array
such that each element has Morley rank 1.
\begin{enumerate}
\item $l$ is a line in $F$ if $l=\{f_{i,j,k}.f_{i',j',k'}\}$ for some
$(i,j,k)\neq(i',j',k')\in\mathbb{N}^{3}$.
\item The unit cube of $F$ is $A=\{f_{i,j,k}\,|\,1\leq i,j,k\leq2\}$.
\item A frame to the unit cube is one of the following sets of lines: \\
\begin{eqnarray*}
 &  & \begin{array}{c}
\mathcal{I}_{1}=\{\{f_{1,1,1},f_{1,1,2}\},\{f_{1,2,1},f_{1,2,2}\},\{f_{2,1,1},f_{2,1,2}\},\{f_{2,2,1}f_{2,2,2}\}\}\\
\mathcal{I}_{2}=\{\{f_{1,1,1},f_{1,2,1}\},\{f_{1,1,2},f_{1,2,2}\},\{f_{2,1,1},f_{2,2,1}\},\{f_{2,1,2}f_{2,2,2}\}\}\\
\mathcal{I}_{3}=\{\{f_{1,1,1},f_{2,1,1}\},\{f_{1,2,1},f_{2,2,1}\},\{f_{1,1,2},f_{2,1,2}\},\{f_{1,2,2}f_{2,2,2}\}\}
\end{array}
\end{eqnarray*}

\item For every frame to the unit cube $\mathcal{I}$ ,we define the base
to the frame to be: $\mathcal{C}_{\mathcal{I}}=\bigcap_{l\in\mathcal{I}}acl(l)$.\\
Let $RM(\mathcal{C}_{j}):=max\{RM(g)\,|\,g\in acl(\mathcal{C}_{j})\}$
 .
\end{enumerate}
\end{defn}
\begin{rem}
\label{Remaek 1.3}If $F=(f_{i,j,k})$ is a 3 indiscernible array
and $l_{1}\neq l_{2}$ are two lines in $F$, then $RM(l_{1},l_{2})>2,$
otherwise $acl(l_{1})=acl(l_{2})$; in contradiction to the fact that
$F$ is a 3 indiscernible array. Hence, if we take three different
lines $l_{1},l_{2},l_{3}$ , with $\mbox{\ensuremath{e\in acl(l_{1})\cap acl(l_{2})\cap acl(l_{3})}}$,
$\mbox{\ensuremath{e'\in acl(l_{1})\cap acl(l_{2})}}$ and $\mbox{\ensuremath{e''\in acl(l_{1})\cap acl(l_{3})}}$,
such that $\mbox{\ensuremath{RM(e)=RM(e')=RM(e'')=1}}$ then we must
have: $\mbox{\ensuremath{acl(e)=acl(e')=acl(e'')}}$.\\

\end{rem}
We can differentiate between the three forms in definition \ref{Definition}
by the Morley rank of the bases to the frames.
\begin{lem}
\label{lem:dividing line}Let $F=(f_{i,j,k})$ be a 3 indiscernible
array.
\begin{enumerate}
\item If $F$ has the group form, then for every frame to the unit cube
$\mathcal{I}$ we have: $RM(\mathcal{C}_{\mathcal{I}})=1$.
\item If $F$ has the field form, then there is exactly one frame to the
unit cube $\mathcal{I}$ such that: $RM(\mathcal{C}_{\mathcal{I}})=1$.
\item If F has the twisted form, then for every frame to the unit cube $\mathcal{I}$
we have: $RM(\mathcal{C}_{\mathcal{I}})=0$.

\begin{proof}
(1)+(2) are trivial; we prove (3):\\
Let $f'_{i,j,k}=\frac{u_{i}+v_{j}}{v_{j}+w_{k}}$, then $\mbox{\ensuremath{\frac{u_{2}+v_{1}}{u_{1}+v_{1}}\in acl(f'_{1,1,1},f'_{2,1,1})\cap acl(f'_{2,1,2},f'_{1,1,2})}}$,
on the other hand: 
\begin{eqnarray*}
 &  & \frac{f'_{2,1,1}-1}{f'_{1,1,1}-1}=\frac{u_{2}-w_{1}}{u_{1}-w_{1}}=\frac{f'_{2,2,1}-1}{f'_{1,2,1}-1}\in acl(f'_{2,1,1},f'_{1,1,1})\cap acl(f'_{2,2,1},f'_{1,2,1}).
\end{eqnarray*}
Suppose there is an element $g\in\mathcal{U}$ such that $\mbox{\ensuremath{RM(g)=1}}$
and $g\in\bigcap_{l\in\mathcal{I}_{3}}acl(l)$, then by remark \ref{Remaek 1.3}
$acl(\frac{u_{2}-w_{1}}{u_{1}-w_{1}})=acl(\frac{u_{2}+v_{1}}{u_{1}+v_{1}})$,
but $acl(\frac{u_{2}-w_{1}}{u_{1}-w_{1}})=acl(\frac{u_{2}-u_{1}}{u_{1}-w_{1}})$
and $acl(\frac{u_{2}+v_{1}}{u_{1}+v_{1}})=acl(\frac{u_{2}-u_{1}}{u_{1}+v_{1}})$
and we have: $w_{1}\in acl(u_{1},u_{2},v_{1})$ in contradiction to
the assumption that $\{w_{1},u_{1},u_{2},v_{1}\}$ is an independent
set.\\
\\
The proof that there is no element $g\in\mathcal{U}$ such that $RM(g)=1$
and $\mbox{\ensuremath{g\in\bigcap_{l\in\mathcal{I}_{i}}acl(l)}}$
for $1\leq i\leq2$ is the same as the proof for $\mathcal{I}_{3}$
if we replace $f'_{i,j,k}$ by $\widehat{f'_{i,j,k}}=(f'_{i,j,k})^{-1}$
for $\mathcal{I}_{1}$ and $\overline{f'_{i,j,k}}=(f'_{i.j.k}-1)^{-1}$
for $\mathcal{I}_{2}$.
\end{proof}
\end{enumerate}
\end{lem}
We conclude this section with the following corollary which answers
Problem \ref{conjecture 1}:
\begin{cor}
Let $\mathbb{F}$ be a definable algebraic closed field, and let $\mbox{\{\ensuremath{a_{i}^{1}\}_{i\in\mathbb{N}}},\{\ensuremath{a_{j}^{2}\}_{j\in\mathbb{N}}},\{\ensuremath{a_{k}^{3}\}_{k\in\mathbb{N}}\subset\mathbb{F}}}$,
where $\{a_{i}^{1},a_{j}^{2},a_{k}^{3}\,|\,i,j,k\in\mathbb{N}\}$
is an algebraic independent set. Define $f_{i,j,k}=\frac{a_{i}^{1}+a_{j}^{2}}{a_{k}^{3}+a_{j}^{2}}$,
then $F=(f_{i,j,k})$ is a counterexample to problem \ref{conjecture 1}.
\begin{proof}
$F$ is a 3 indiscernible array (because of the indiscerniblity of
the set $\{a_{i}^{1},a_{j}^{2},a_{k}^{3}\,|\,i,j,k\in\mathbb{N}\}$).
By lemma \ref{lem: equal 2}, $\alpha_{F}(m,n,p)=m+n+p-2$, and by
lemma \ref{lem:dividing line},  $F$ neither has the group form nor
the field form.
\end{proof}
\end{cor}

\section{Indiscernible arrays generated by definable fields}

In this section we prove some basic facts about indiscernible arrays
generated by fields or groups.
\begin{lem}
\label{lem:field}Let $(\mathbb{F},\cdot,+)$ be a definable field,
and $F=(f_{i,j})$ an indiscernible array such that $\mbox{\ensuremath{f_{i,j}}=\ensuremath{a_{i}\cdot c_{j}}+\ensuremath{b_{j}}}$,
where $b_{j}\in acl(c_{j})$, $b_{1}=0$ and \textup{$c_{1}=1$,}
then there exists an element $k\in acl(\emptyset)\cap\mathbb{F}$
such that $f_{i,j}=a_{i}c_{j}+k(c_{j}-1)$. 
\begin{proof}
Let $v_{k,l}=\frac{c_{l}}{c_{k}}$ and $t_{k,l}=c_{k}-v_{k.l}\cdot c_{l}$.
If $t_{k,l}=0$ for some $k,l\in\mathbb{N}$, then by indiscerniblity
it is true for every $k,l$, thus: $f_{i,j}=a_{i}c_{j}$. Hence, we
may assume $t_{k,l}\neq0$. On the other hand: $v_{k,l}=\frac{f_{1,l}-f_{2,l}}{f_{1,k}-f_{2,k}}$
and $t_{k,l}=f_{1,k}-f_{1,l}\cdot v_{k,l}$, thus: $(1)$ $v_{k,l},t_{k,l}\in dcl(f_{i,j}\,|\,\ j\in\{l,k\},\quad1\leq i\leq2).$\\

\begin{claim*}
1: for every $k<l\in\mathbb{N}$: $t_{k,l}\in acl(v_{k,l})$.
\begin{proof}
By our assumptions $v_{1,2}=b_{2},\quad t_{1,2}=c_{2}$ and $b_{2}\in acl(c_{2})$.
From indiscernibility 
\begin{eqnarray*}
\  &  & tp(f_{i,j}\,|\,\ j\in\{l,k\},\,\,1\leq i\leq2)=tp(f_{i,j}\,|\,1\leq j\leq2,\,\,1\leq i\leq2),
\end{eqnarray*}
 thus from $(1)$ we get $t_{k,l}\in acl(v_{k,l})$.
\end{proof}
\end{claim*}
Let $N\in\mathbb{N}$ be the number of conjugates of $t_{1,2}$ over
$v_{1,2}$ and fix $n>N+2$. We  look at 
\begin{eqnarray*}
 &  & ((v_{1,2},t_{1,2}),(v_{2,3},t_{2,3}),..,(v_{n,n+1},t_{n,n+1})).
\end{eqnarray*}

\begin{claim*}
2: For every permutation $\sigma\in S_{n}$ , there is an automorphism
$\tau_{\sigma}$ which sends:\\
 $\mbox{V=((\ensuremath{v_{1,2}},\ensuremath{t_{1,2}}),(\ensuremath{v_{2,3}},\ensuremath{t_{2,3}}),..,(\ensuremath{v_{n,n+1}},\ensuremath{t_{n,n+1}}))}$
to $\mbox{\ensuremath{V_{\sigma}}=((\ensuremath{v_{\sigma(1),\sigma(1)+1}},\ensuremath{t_{\sigma(1),\sigma(1)+1}}),..,(\ensuremath{v_{\sigma(n),\sigma(n+1)}},\ensuremath{t_{1,n}}))}$.
\begin{proof}
Let $p=tp(v_{1,2},t_{1,2})$, and define 
\begin{eqnarray*}
 &  & A=\{(\overline{s_{1}},\overline{s_{2}},...,\overline{s_{n}})\,|\,\overline{s_{i}}\vDash p,\,\,\{\overline{s_{1}},..,\overline{s_{n}}\}-is\,\,independent\,\,set\}\ 
\end{eqnarray*}
 by uniqueness of non-forking extensions of types, there is a complete
type $q(\overline{x_{1}},..,\overline{x_{n}})$ such that: $A=\{(\overline{s_{1}},\overline{s_{2}},...,\overline{s_{n}})\,|\,(\overline{s_{1}},\overline{s_{2}},...,\overline{s_{n}})\vDash q\}.$
$V,V_{\sigma}\in A$, thus there exists an automorphism $\tau_{\sigma}$
such that $\tau_{\sigma}(V)=V_{\sigma}$. \end{proof}
\begin{claim*}
3 Let $\sigma\in S_{n}$,\end{claim*}
\begin{enumerate}
\item If there are $1\leq M_{1}<M_{2}\leq n$ such that $\sigma(i)=i$ for
every $M_{1}\leq i\leq M_{2}$ , then $\mbox{\ensuremath{\tau_{\sigma}}(\ensuremath{t_{M_{1},M_{2}}})=\ensuremath{t_{M_{1},M_{2}}}}$
and $\mbox{\ensuremath{\tau_{\sigma}}(\ensuremath{v_{M_{1},M_{2}}})=\ensuremath{v_{M_{1},M_{2}}}}$.
\item If there exists an $M<n$ such that for every $M\leq i\leq n$ we
have: $\sigma(i)=i$, and $\tau_{\sigma}(t_{1,n+1})=t_{1,n+1}$, then
$\tau_{\sigma}(t_{1,M})=t_{1,M}$.\end{enumerate}
\begin{proof}
$ $
\begin{enumerate}
\item Observe that $\mbox{\ensuremath{\tau_{\sigma}}(\ensuremath{v_{M_{1},M_{2}}})=\ensuremath{\tau_{\sigma}}(\ensuremath{\prod_{i=M_{1}}^{M_{2}-1}v_{i.i+1}})=\ensuremath{\prod_{i=M_{1}}^{M_{2}-1}\tau_{\sigma}}(\ensuremath{v_{i.i+1}})=\ensuremath{\prod_{i=M_{1}}^{M_{2}-1}v_{i.i+1}}=\ensuremath{v_{M_{1},M_{2}}}}$.
We prove that $\tau_{\sigma}(t_{M_{1},M_{2}})=t_{M_{1},M_{2}}$ by
induction on $M_{2}-M_{1}$. The basis of the induction $(M_{2}-M_{1}=1)$
follows by our assumption that $\sigma(M_{1})=M_{1}$, for the induction
step:
\begin{eqnarray*}
 &  & \tau_{\sigma}(t_{M_{1},M_{2}})=\tau_{\sigma}(t_{M_{1},M_{2}-1})\tau_{\sigma}(v_{M_{2}-1,M_{2}})+\tau_{\sigma}(t_{M_{2}-1,M_{2}})=t_{M_{1},M_{2}-1}v_{M_{2}-1,M_{2}}+t_{M_{2}-1,M_{2}}=t_{M_{1},M_{2}}
\end{eqnarray*}
 (by induction hypothesis $\tau_{\sigma}(t_{M_{1},M_{2}-1})=t_{M_{1},M_{2}-1}$).
\item We have $t_{1,n+1}=t_{1,M}v_{M,n+1}+t_{M,n+1}$, by the first part
of the claim $\tau_{\sigma}(t_{M,n+1})=t_{M,n+1}$ and $\mbox{\ensuremath{\tau_{\sigma}}(\ensuremath{v_{M_{1},M_{2}}})=\ensuremath{v_{M_{1},M_{2}}}}$
thus $\tau_{\sigma}$ fix $t_{1,M}$.
\end{enumerate}
\end{proof}
\end{claim*}

We are now ready to prove the main claim:
\begin{claim*}
4: For every $k\in\mathbb{N}$: $t_{k,k+2}=t_{k,k+1}v_{k+1,k+2}+t_{k+1,k+2}=t_{k+1,k+2}v_{k,k+1}+t_{k,k+1}$.
\begin{proof}
Let $\sigma_{i}=(i-1,i)\in S_{n}$, that is the permutation which
replaces $i-1$ and $i$. By claim 1, $\tau_{\sigma_{l}}(t_{1,n+1})\in acl(v_{1,n+1})$
and because $n>N$ there must exists $1<l<m\leq n$ such that $\mbox{\ensuremath{\tau_{\sigma_{l}}}(\ensuremath{t_{1,n+1}})=\ensuremath{\tau_{\sigma_{m}}}(\ensuremath{t_{1,n+1}})}$.
Hence, by Claim $3(2)$: $\tau_{\sigma_{l}}(t_{1,m+1})=\tau_{\sigma_{m}}(t_{1,m+1})$
and using Claim 3 once again give us:
\begin{eqnarray*}
(4) &  & \tau_{\sigma_{l}}(t_{l-1,l+1})-t_{l-1,l+1}=\tau_{\sigma_{l}}(t_{1,l+1})-t_{1,l+1}=\frac{\tau_{\sigma_{l}}(t_{1,m+1})-t_{1,m+1}}{v_{l+1,m+1}}=\frac{\tau_{\sigma_{m}}(t_{1,m+1})-t_{1,m+1}}{v_{l+1,m+1}}=\\
 &  & \frac{\tau_{\sigma_{m}}(t_{l+1,m+1})-t_{l+1,m+1}}{v_{l+1,m+1}}=\frac{\tau_{\sigma_{m}}(t_{m-1,m+1})-t_{m-1,m+1}}{v_{l+1,m+1}}=\alpha.
\end{eqnarray*}
 However, by (3) $\frac{\tau_{\sigma_{m}}(t_{l+1,m+1})-t_{l+1,m+1}}{v_{l+1,m+1}}\in acl(v_{l+1,m+1})$
and $\tau_{\sigma_{l}}(t_{1,l+1})-t_{1,l+1}\in acl(v_{1,l+1})$. We
also know that $v_{1,l+1}\downarrow v_{l+1,m+1}$ (or else $\{f_{j,m+1}\}_{j\in\mathbb{N}}\in acl\{f_{i,1},f_{i,l+1}\,|\,i\in\mathbb{N}\}$
in contradiction to the fact that $F$ is an indiscernible array),
therefore we can conclude that $\alpha\in acl(\emptyset)$.\\
\\
On the one hand, by equation$(4)$: 
\begin{eqnarray*}
 &  & t_{l,l+1}v_{l-1,l}+t_{l-1,l}=\tau_{\sigma_{l}}(t_{l-1,l+1})=\alpha+t_{l-1,l+1}=\alpha+t_{l-1,l}v_{l,l+1}+t_{l,l+1}.
\end{eqnarray*}
On the other hand: 
\begin{eqnarray*}
 &  & t_{m,m+1}v_{m-1,m}+t_{m-1,m}=\tau_{\sigma_{m}}(t_{m-1,m+1})=\alpha v_{l+1,m+1}+t_{m-1,m+1}=\\
 &  & \alpha\cdot v_{l+1,m+1}+t_{m-1,m}v_{m,m+1}+t_{m,m+1}.
\end{eqnarray*}
However, $tp((t_{l-1,l},v_{l,l-1}),(t_{l,l+1},v_{l.l+1}))=tp((t_{m-1,m},v_{m,m-1}),(t_{m,m+1},v_{m,m+1}))$
thus necessarily:
\begin{eqnarray*}
 &  & \alpha\cdot v_{l+1,m+1}+t_{m-1,m}v_{m,m+1}+t_{m,m+1}=\alpha+t_{m-1,m}v_{m,m+1}+t_{m,m+1},
\end{eqnarray*}
 so $\alpha=0$, and we get the desired result: $t_{l,l+1}v_{l-1,l}+t_{l-1,l}=t_{l-1,l}v_{l,l+1}+t_{l,l+1}$.
\end{proof}

$ $\\
Finally we can prove our lemma: By claim 4 $\frac{t_{k+1,k+2}}{(v_{k+1,k+2}-1)}=\frac{t_{k,k+1}}{(v_{k,k+1}-1)}=\alpha$.
Therefore, 
\begin{eqnarray*}
 &  & \alpha\in acl(v_{k+1,k+2})\cap acl(v_{k,k+1})\cap\mathbb{F}=acl(\emptyset)\cap\mathbb{F}
\end{eqnarray*}
 as desired.

\end{claim*}
\end{proof}

\newpage{}
\begin{lem}
\textup{\label{lem:Field2}Let $(\mathbb{F},\cdot,+)$ be a definable
field.}\end{lem}
\begin{enumerate}
\item If $F=(f_{i,j})$ \textup{is an indiscernible array, and there are
elements $\mbox{\{\ensuremath{u_{i}\}_{i\in\mathbb{N}}},\{\ensuremath{v_{j}\}_{j\in\mathbb{N}}},\{\ensuremath{w_{k}\}_{k\in\mathbb{N}}\subset\mathbb{F}}},$
such that: $v_{j}\in acl(\emptyset),\quad v_{1}=1$ and $f_{i,j}=u_{i}v_{j}+w_{j}$,
then $v_{j}=1$ for every $i,j\in\mathbb{N}$.}
\item If $F=(f_{i,j,k})$ \textup{i}s a \textup{3 indiscernible array, and
there are} \textup{elements $\mbox{\{\ensuremath{u_{i}\}_{i\in\mathbb{N}}},\{\ensuremath{v_{j,k}\}_{j,k\in\mathbb{N}}},\{\ensuremath{w_{j,k}\}_{j,k\in\mathbb{N}}\subset\mathbb{F}}}$
}\\
\textup{such that for every $i,j,k\in\mathbb{N}$ we have: $\mbox{\ensuremath{v_{1,k}}=1,\quad\ensuremath{w_{1,1}}=0,\quad\ensuremath{w_{j,1}\in}acl(\ensuremath{v_{j,1}})}$
and $\mbox{\ensuremath{f_{i,j,k}}=\ensuremath{u_{i}v_{j,k}}+\ensuremath{w_{j,k}}},$
then $F$ has the field form.}
\item If $F=(f_{i,j,k})$ \textup{is a} \textup{3 indiscernible array and
there are elements $\mbox{\{\ensuremath{u_{i}\}_{i\in\mathbb{N}}},\{\ensuremath{v_{j,k}\}_{j,k\in\mathbb{N}}},\{\ensuremath{w_{j,k}\}_{j,k\in\mathbb{N}}\subset\mathbb{F}}}$
}\\
\textup{such that for every $i,j,k\in\mathbb{N}$ we have: $\mbox{\ensuremath{v_{1,1}}=1,\quad\ensuremath{w_{1,1}}=0,\quad\ensuremath{w_{j,1}\in}acl(\ensuremath{v_{j,1}}),\quad\ensuremath{w_{1,k}\in}acl(\ensuremath{v_{1,k}})}$,
$f_{i,j,k}=u_{i}v_{j,k}+w_{j,k}$, and }$RM(\mathcal{C}_{\mathcal{I}_{3}})=0$
\textup{where:
\begin{eqnarray*}
 &  & \mathcal{I}_{3}=\{\{f_{1,1,1},f_{2,1,1}\},\{f_{1,2,1},f_{2,2,1}\},\{f_{1,1,2},f_{2,1,2}\},\{f_{1,2,2}f_{2,2,2}\}\}
\end{eqnarray*}
 (see Definition }\ref{def:frames}\textup{)}, \textup{then $F$ has
the field form, or $F$ has the twisted form.} \end{enumerate}
\begin{proof}
$ $
\begin{enumerate}
\item Let $\sigma_{3}$ be an automorphism such that $\sigma_{3}(f_{i,1})=f_{i,1}$
and $\sigma_{3}(f_{i,2})=f_{i,3}$ for every $1\leq i\leq2$, then:
$v_{2}=\sigma_{3}(v_{2})=v_{3}$.\\
 \\
For every $t\in\mathbb{N}$ let $\sigma_{t}'$ be an automorphism
such that: $\sigma'_{t}(f_{i,2})=f_{i,1}$ and $\sigma_{t}'(f_{i,3})=f_{i,t}$
for every $1\leq i\leq2$, then: 
\begin{eqnarray*}
 &  & 1=\sigma'_{t}(1)=\sigma'_{t}(\frac{v_{3}}{v_{2}})=\sigma'_{t}(\frac{f_{2,3}-f_{1,3}}{f_{2,2}-f_{1,2}})=\frac{f_{2,t}-f_{1,t}}{f_{2,1}-f_{1,1}}=\frac{v_{t}}{v_{1}}=v_{t}.
\end{eqnarray*}

\item By Lemma \ref{lem:field} there exist $\alpha\in acl(\emptyset)\cap\mathbb{F}$
such that $f_{i,j,1}=u{}_{i}v_{j,1}+\alpha(v{}_{j,1}-1)$, thus by
\\
indiscerniblity: $\frac{f_{2,1,k}+\alpha}{f_{1,1,k}+\alpha}=\frac{f_{2,j,k}+\alpha}{f_{1,j,k}+\alpha}=\frac{u_{2}+\alpha}{u_{1}+\alpha}$
for every $k\in\mathbb{N}$. By our assumption 
\begin{eqnarray*}
 &  & f_{i,j,k}-f_{i,j,1}=f_{1,j,k}-f_{1,j,1}=0
\end{eqnarray*}
 for $j=1$, and again by indiscerniblity, this is true for every
$j\in\mathbb{N}$. Let $\mbox{\ensuremath{f_{1,j,k}}-\ensuremath{f_{1,1,1}}=\ensuremath{t_{k}^{j}}}$,
notice that:
\begin{eqnarray*}
 &  & \mbox{\ensuremath{\frac{f_{2,1,1}+t_{k}^{1}+\alpha}{f_{1,1,1}+t_{k}^{1}+\alpha}}=\ensuremath{\frac{f_{2,1,k}+\alpha}{f_{1,1,k}+\alpha}}=\ensuremath{\frac{f_{2,j,k}+\alpha}{f_{1,j,k}+\alpha}}=\ensuremath{\frac{f_{2,j,1}+\alpha+t_{k}^{j}}{f_{1,j,1}+\alpha+t_{k}^{j}}}.}
\end{eqnarray*}
 However, $f_{i,j,1}+\alpha=(f_{i,1,1}+\alpha)\cdot v_{j,1}$, thus:
\begin{eqnarray*}
 &  & \frac{f_{2,1,1}+t_{k}^{1}+\alpha}{f_{1,1,1}+t_{k}^{1}+\alpha}=\frac{(f_{2,1,1}+\alpha)\cdot v_{j,1}+t_{k}^{j}}{(f_{1,1,1}+\alpha)\cdot v_{j,1}+t_{k}^{j}},
\end{eqnarray*}
which implies: $t_{k}^{j}=v_{j,1}t_{k}^{1}$ . \\
\\
Let $f'_{i,j,k}=f_{i,j,k}+\alpha,\,\,u'_{i}=u_{i}+\alpha,\,\,v'_{j}=v_{j,1},\,\,w'_{k}=w_{1,k}=t_{k}^{1}$,
then: $f'_{i,j,k}\in acl(f_{i,j,k})$ and we have:
\begin{eqnarray*}
 &  & f'_{i,j,k}=f_{i,j,k}+\alpha=f_{i,j,1}+t_{k}^{j}+\alpha=f_{i,j,1}+t_{k}^{1}v_{j,1}+\alpha=(u_{i}+\alpha+t_{k}^{1})v_{j,1}=(u'_{i}+w'_{k})\cdot v'_{j}.
\end{eqnarray*}

\item By Lemma \ref{lem:field} there are $\alpha,\beta\in acl(\emptyset)\cap\mathbb{F}$
such that $f_{i,j,1}=u{}_{i}v{}_{j,1}+\alpha(v{}_{j,1}-1)$ and \\
$f_{i,1,k}=u{}_{i}v{}_{j,1}+\beta(v{}_{1,k}-1)$ $(\dagger)$. \\
\\
If $\alpha=\beta$, then: 
\begin{eqnarray*}
 &  & \mbox{\ensuremath{\frac{f_{2,2,2}+\alpha}{f_{1,2,2}+\alpha}}=\ensuremath{\frac{f_{2,1,1}+\alpha}{f_{1,1,1}+\alpha}}=\ensuremath{\frac{f_{2,1,2}+\alpha}{f_{1,1,2}+\alpha}}=\ensuremath{\frac{f_{2,2,1}+\alpha}{f_{1,2,1}+\alpha}}=\ensuremath{\frac{u_{2}+\alpha}{u_{1}+\alpha}}}
\end{eqnarray*}
and of course $RM(\frac{f_{2,1,1}+\alpha}{f_{1,1,1}+\alpha})=1$,
thus $RM(\mathcal{C}_{\mathcal{I}_{3}})=1$, in contradiction to our
assumption.\\
 \\
Assume $\alpha\neq\beta$: Let $\widetilde{u_{i}^{j}}=f_{i,j,1},\,\,v_{k}^{j}=\frac{v_{j,k}}{v_{j,1}}$
and $\widetilde{w_{k}^{j}}=w_{j,k}-w_{j,1}\cdot\widetilde{v_{k}^{j}}$,
then $\mbox{\ensuremath{f_{i,j,k}}=\ensuremath{\widetilde{u_{i}^{j}}\cdot\widetilde{v_{k}^{j}}}+\ensuremath{\widetilde{w_{k}^{j}}}.}$
\\
Observe that for $j=1$ we have: $\mbox{\ensuremath{\widetilde{w_{k}^{1}}}=\ensuremath{w_{1,k}}=\ensuremath{\beta}(\ensuremath{v_{1,k}}-1)=\ensuremath{\beta}(\ensuremath{\widetilde{v_{k}^{1}}}-1).}$
\\
\\
Let $\sigma_{j}$ be an automorphism which takes $\mbox{\ensuremath{(f_{i,1,l}\,|\,1\leq i\leq2,\quad l\in\{1,k\})}}$
to $\mbox{\ensuremath{(f_{i,j,l}\,|\,1\leq i\leq2,\quad l\in\{1,k\})}}$,
then it satisfies: $\widetilde{w_{k}^{j}}=\sigma_{j}(\widetilde{w_{k}^{1}})=\beta(\widetilde{v_{k}^{j}}-1)$.
In conclusion we have: $(1)$$f_{i,j,k}=f_{i,j,1}\frac{v_{j,k}}{v_{j,1}}+\frac{v_{j,k}}{v_{j,1}}\beta-\beta$,
and in the same way we get: $(2)$$f_{i,j,k}=f_{i,1,k}\frac{v_{j,k}}{v_{1,k}}+\frac{v_{j,k}}{v_{1,k}}\alpha-\alpha$.
\\
\\
We can rewrite $(1)$ as follows: 
\begin{eqnarray*}
 &  & f_{1,j,k}=f_{1,j,1}\frac{v_{j,k}}{v_{j,1}}+\frac{v_{j,k}}{v_{j,1}}\beta-\beta=(f_{1,1,1}v_{j,1}+v_{j,1}\alpha-\alpha)\frac{v_{j,k}}{v_{j,1}}+\frac{v_{j,k}}{v_{j,1}}\beta-\beta.
\end{eqnarray*}
On the other hand, by $(2)$ we get: 
\begin{eqnarray*}
 &  & f_{1,j,k}=f_{1,1,k}\frac{v_{j,k}}{v_{1,k}}+\frac{v_{j,k}}{v_{1,k}}\alpha-\alpha=(f_{1,1,1}v_{1,k}+v_{1,k}\beta-\beta)\frac{v_{j,k}}{v_{1,k}}+\frac{v_{j,k}}{v_{1,k}}\alpha-\alpha,
\end{eqnarray*}
thus: $(\alpha-\beta)(v_{j,k}-\frac{v_{j,k}}{v_{j,1}}-\frac{v_{j,k}}{v_{1,k}}+1)=0$
. Therefore, by our assumption:
\begin{eqnarray*}
 &  & 0=\frac{v_{j,k}}{v_{j,1}}\cdot v_{j,1}-\frac{v_{j,k}}{v_{j,1}}-\frac{v_{j,k}}{v_{1,k}}+1=\frac{v_{j,k}}{v_{j,1}}\cdot v_{j,1}-\frac{v_{j,k}}{v_{j,1}}-\frac{\left(\frac{v_{j.k}}{v_{j,1}}\right)}{v_{1,k}}v_{j,1}+1,
\end{eqnarray*}
which implies $\frac{v_{j,k}}{v_{j,1}}=\frac{v_{1,k}}{v_{1,k}-v_{1,k}v_{j,1}+v_{j,1}}$
and 
\begin{eqnarray*}
 &  & f_{i,j,k}=(f_{i,1,1}v_{j,1}+v_{j,1}\alpha-\alpha)\frac{v_{1,k}}{v_{1,k}-v_{1,k}v_{j,1}+v_{j,1}}+\frac{v_{1,k}}{v_{1,k}-v_{1,k}v_{j,1}+v_{j,1}}\beta-\beta.
\end{eqnarray*}
Let: $f'_{i,j,k}=\frac{f_{i,j,k}+\alpha}{\beta-\alpha},\,\,w'_{k}=v_{1,k},\,\,v'_{j}=v_{j,1},\,\,u'_{i}=f'_{i,1,1}$
, then: $acl(f'_{i,j,k})=acl(f_{i,j,k})$ and 
\begin{eqnarray*}
 & f'_{i,j,k} & =\frac{(f_{i,1,1}+\alpha)v'_{j}-\alpha+\beta)\frac{w'_{k}}{w'_{k}-w'_{k}v'_{j}+v'_{j}}+\alpha-\beta}{\beta-\alpha}=(\frac{f_{i,1,1}+\alpha}{\beta-\alpha})v'_{j}+1)\frac{w'_{k}}{w'_{k}-w'_{k}v'_{j}+v'_{j}}-1\\
 &  & =(u'_{i}v'_{j}+1)\frac{w'_{k}}{w'_{k}-w'_{k}v'_{j}+v'_{j}}-1
\end{eqnarray*}
Let $v''_{j}=\frac{1}{v'_{j}},\quad w''_{k}=\frac{1}{w'_{k}}-1$ and
$t_{i,j,k}=(u'_{i}v'_{j}+1)\frac{w'_{k}}{w'_{k}-w'_{k}v'_{j}+v'_{j}}$,then
$t_{i,j,k}=\frac{u'_{i}+v''_{j}}{v''_{j}+w''_{k}}$ and \\
$acl(t_{i,j,k})=acl(f_{i,j,k})$. Thus $F$ has the twisted form. 
\end{enumerate}
\end{proof}
\end{lem}

\section{Main Theorem}

In the previous section we studied indiscernible arrays that were
induced by fields/group, the following lemma will be helpful in retrieving
indiscernible arrays of this form. 
\begin{lem}
\label{lem:Constructing lemma} Let $F=(f_{i,j,k})$ be a 3 indiscernible
array, $\{a_{i}\}_{i\in\mathbb{N}}\subset\mathcal{U},\quad\{\overline{b_{j,k}}\}_{j,k\in\mathbb{N}}\subset\mathcal{U}^{m}$
and $h:\,\mathcal{U}^{m+1}\rightarrow\mathcal{U}$ a definable function.
\begin{enumerate}
\item If there is a formula $\phi(x,y)$ such that for every $i,j,k\in\mathbb{N}$:
$h(a_{i},\overline{b}_{j,k})\vDash\phi(x,f_{i,j,k})$ and\\
$|\{a\vDash\phi(x,f_{i,j,k})\,|\,\ a\in\mathcal{U}\}|<\infty$, then
there are $\mbox{\{a\ensuremath{'_{i}\}_{i\in\mathbb{N}}\subset\mathcal{U}},\quad\{\ensuremath{\overline{b'_{j,k}}\}_{j,k\in\mathbb{N}}\subset\mathcal{U}^{m}}}$
and a 3 \\
indiscernible array $F'=(f'_{i,j,k})$ such that $f'_{i,j,k}=h(a'_{i},\overline{b'_{j,k}})$
and $F'$ is algebraically equivalent \\
to $F$.
\item Suppose there is a formula $\phi(x,y)$ , such that for every $i,j\in\mathbb{N}$:
$\mbox{h(\ensuremath{a_{i}},\ensuremath{\overline{b}_{j,j}})\ensuremath{\vDash\phi}(x,\ensuremath{f_{i,j,j}})}$,\\
 $|\{a\vDash\phi(x,f_{i,j,j})\,|\,\ a\in\mathcal{U}\}|<\infty$ and
$h(a_{i},\overline{b_{1,1}})=a_{i}$, then there are $\mbox{\ensuremath{\{a'_{i}\}_{i\in\mathbb{N}}\subset\mathcal{U},\,\,\{\overline{b'_{j,k}}\}_{j,k\in\mathbb{N}}\subset\mathcal{U}^{m}}}$
and a 3 indiscernible array $F'=(f'_{i,j,k})$ such that: $f'_{i,j,k}=h(a'_{i},\overline{b'_{j,k}})$
and $F'$ is algebraically equivalent to $F$.\end{enumerate}
\begin{proof}
$ $
\begin{enumerate}
\item Let $h(a_{i},\overline{b_{j,k}})=t_{i,j,k}$. We  use the infinite
Ramsey theorem: Fix $N\in\mathbb{N}$, for every $r_{1},..,r_{N}\in\mathbb{N}$
let $\sigma_{r_{1},..,r_{N}}$ be an automorphism such that: $\mbox{\ensuremath{(\sigma(f_{i,j,r_{k}})\,|\,1\leq i,k,j\leq N)=(f_{i,j,k}\,|\,1\leq i,j,k\leq N)}}$.
By the infinite Ramsey theorem there is an infinite subset $I\subset\mathbb{N}$
such that for every $\mbox{\ensuremath{r_{1},..,r_{N},t_{1},..,t_{N}\in I}}$
we have: $\mbox{\ensuremath{tp(t'_{i,j,r_{k}}\,|\,1\leq i,k,j\leq N)=tp(t'_{i,j,t_{k}}\,|\,1\leq i,j,k\leq N)}}$.
By compactness we can find $\{a'_{i}\}_{i\in\mathbb{N}}\subset\mathcal{U},\quad\{\overline{b'_{j,k}}\}_{j,k\in\mathbb{N}}\subset\mathcal{U}^{m}$
and an array $F'=(t'_{i,j,k})$ such that: $\mbox{\ensuremath{t'_{i,j,k}=h(a'_{i},\overline{b'_{j,k}})\in acl(f_{i,j,k})}}$
and for every $n\in\mathbb{N}$ and $r_{1},..,r_{n}\in\mathbb{N}$
we have:
\begin{eqnarray*}
 &  & tp(t'_{i,j,r_{k}}\,|\,1\leq i,k,j\leq n)=tp(t'_{i,j,k}\,|\,1\leq i,j,k\leq n)(\dagger).
\end{eqnarray*}
If we repeat this construction, this time with respect to $(t'_{i,j,k})_{i,j,k\in\mathbb{N}}$
, we can find $\{a''_{i}\}_{i\in\mathbb{N}}\subset\mathcal{U}$ and
$\{\overline{b''_{j,k}}\}_{j,k\in\mathbb{N}}\subset\mathcal{U}^{m}$
such that $\mbox{\ensuremath{t'_{i,j,k}=h(a'_{i},\overline{b'_{j,k}})\in acl(f_{i,j,k})}}$
satisfies $(\dagger)$, and also for every $n\in\mathbb{N}$, and
$r_{1},..,r_{n}\in\mathbb{N}$ we have: $\mbox{\ensuremath{tp(t''_{i,r_{j},k}\,|\,1\leq i,j,k\leq n)=tp(t''_{i,j,k}\,|\,1\leq i,j,k\leq n).}}$
\\
\\
Repeating this construction for the last time with respect to $(t''_{i,j,k})$,
we can find: $\{a^{(3)}\}_{i\in\mathbb{N}}\subset\mathcal{U},$ $\{\overline{b^{(3)}{}_{j,k}}\}_{j,k\in\mathbb{N}}\subset\mathcal{U}^{m}$
and an array $F''=(f'_{i,j,k})$ (where $f'_{i,j,k}=h(a_{i}^{(3)},\overline{b_{j,k}^{(3)}})$)
which satisfies $(\dagger)$ and also for every $n\in\mathbb{N}$
and for every $r_{1},..,r_{n}\in\mathbb{N}$ : 
\begin{eqnarray*}
 &  & tp(f'_{i,r_{j},k}\,|\,1\leq i,j,k\leq n)=tp(f'_{i,j,k}\,|\,1\leq i,j,k\leq n)\\
 &  & tp(f'_{r_{i},j,k}\,|\,1\leq i,j,k\leq n)=tp(f'_{i,j,k}\,|\,1\leq i,j,k\leq n)
\end{eqnarray*}
Thus for every $n\in\mathbb{N}$ and for every $r_{1},..,r_{n},t_{1},..,t_{n},w_{1},..,w_{n}\in\mathbb{N}$:
\begin{eqnarray*}
 &  & tp(f'_{r_{i},t_{j},w_{k}}\,|\,1\leq i,j,k\leq n)=tp(f'_{i,t_{j},w_{k}}\,|\,1\leq i,j,k\leq n)=\\
 &  & tp(f'_{i,j,w_{k}}\,|\,1\leq i,j,k\leq n)=tp(f'_{i,j,k}\,|\,1\leq i,j,k\leq n)
\end{eqnarray*}
therefore, $F'$ is a 3 indiscernible array.
\item Let $M=|\{a\vDash\phi(x,f_{i,j,k})\,|\,a\in\mathcal{U}\}|$. Fix $N\in\mathbb{N}$,
by indecerniblity for every $\mbox{\ensuremath{1\leq i\leq N,\,\,1\leq j\leq(N+1)\cdot M^{N}}}$
and $k>(N+1)\cdot M^{N}$ there is an automorphism $\sigma_{j,k}$
such that: $\mbox{\ensuremath{\sigma_{j,k}}(\ensuremath{f_{i,1,1}})=\ensuremath{f_{i,1,1}}}$
and $\mbox{\ensuremath{\sigma_{j,k}}(\ensuremath{f_{i,j,j}})=\ensuremath{f_{i,j,k}}}$.
By our assumptions for every $1\leq i\leq N$ and $1\leq j\leq M^{N}+1$:
$|\{\sigma_{j,k}(a_{i})\,|k>M^{N}+1\}|<\infty$, thus there is an
infinite subset $\mbox{\ensuremath{I_{1}\subset\mathbb{N}\backslash}\{1,..,(N+1)\ensuremath{\cdot M^{N}}\}}$
such that for every $1\leq i\leq N$ and for every $r,l\in I_{1}$:
$\sigma_{2,r}(a_{i})=\sigma_{2,l}(a_{i})$.\\
\\
We resume and build $I_{j}$ constructively such that $I_{j}\subset I_{j-1}$
is an infinite subset and $\mbox{\ensuremath{\sigma_{j,r}}(\ensuremath{a_{i}})=\ensuremath{\sigma_{j,l}}(\ensuremath{a_{i}})}$
for every $1\leq i\leq N$ and every $r,l\in I_{j}$. By the pigeonhole
principle there is a set $\mbox{\ensuremath{J\subset\{1,..,(N+1)\cdot M^{N}\}}}$
such that $|J|\geq N$ and for every $1\leq i\leq N$, $j,j'\in J$
and $r,r'\in I_{M^{n}+1}$ : $\sigma_{j,r}(a_{i})=\sigma_{j',r'}(a_{i})=a'_{i}.$
In particular $\mbox{\ensuremath{h(a'_{i},\sigma_{j,k}(\overline{b_{j,j}}))\in acl(f_{i,j,k})}}$
for every $1\leq i\leq N,\quad j\in J,\quad k\in I_{(N+1)\cdot M^{N}}$.
By \\
compactness we can find $\mbox{\ensuremath{\{a''_{i}\}_{i\in\mathbb{N}},\{\overline{b'_{j,k}}\}\in\mathbb{N}}}$
such that: $h(a''_{i}\overline{b'_{j,k}}))\in acl(f_{i,j,k})$. We
may now use the first part of this lemma in order to get the desired
3 indiscernible array.
\end{enumerate}
\end{proof}
\end{lem}
The next lemma takes care of the group form case:
\begin{lem}
\label{lem:group form} Let $F=(f_{i,j,k})$ be a 3 indiscernible
array with\\
 $\alpha(m,n,p)=m+n+p-2$ such that for every frame of the unit cube
$\mathcal{I}$ there is an element $g\in\mathcal{U}$ such that $RM(g)=1$
and $g\in\bigcap_{l\in\mathcal{I}}acl(l)$. Then $F$ has the group
form.
\begin{proof}
Conssider the indiscernible array $F'=(f_{i,j,j}\,|\,i,j\in\mathbb{N})$.
$\mbox{\ensuremath{RM(f_{1,1,1},f_{1,2,2},f_{2,1,1},f_{2,2,2})=3}}$,
thus by \ref{lem:Counter1} (2): $\alpha_{F'}(m,n)=m+n-1$. By fact
\ref{fct:2.2}(1) there is a definable one dimensional group $(G,\cdot)$
and elements $\{b_{i}\}_{i\in\mathbb{N}},\{a_{j}\}_{j\in\mathbb{N}}\subset G$
such that $acl(f_{i,j,j})=acl(a_{i}\cdot b_{j})$ for every $i,j\in\mathbb{N}$.
By Lemma \ref{lem:Constructing lemma} there exist an indiscernible
array $F'=(f'_{i,j,k})$ which is algebraically equivalent to $F$
such that $f'_{i.j.k}=a'_{i}\cdot b'_{j,k}$, where $\mbox{\ensuremath{\{a'_{i}\}_{i\in\mathbb{N}},\{b'_{j,k}\}_{j,k\in\mathbb{N}}\subset G}}$
.\\
For every $i,j,k,j,n\in\mathbb{N}$: $\frac{f'_{i,j',k}}{f'_{i,j,k}}=\frac{b'_{l,k}}{b'_{j,k}}=\frac{f'_{i+1,j',k}}{f'_{i+1,j,k}}$,
and by Remark \ref{Remaek 1.3} $acl(\frac{b'_{j',k}}{b'_{j,k}})=acl(\frac{b'_{j',k'}}{b'_{j,k'}})$
for every $k,k'\in\mathbb{N}$, so for every $j',j\in\mathbb{N}$
there must exists $k,k'\in\mathbb{N}$ such that $\frac{b'_{l,k}}{b'_{j,k}}=\frac{b'_{l,k'}}{b'_{j,k'}}$,
thus by indiscerniblity $\frac{b'_{j',k}}{b'_{j,k}}$ does not depend
on $k$. Fixing $j$ instead of $k$ we get: $\frac{f'_{i,j,k'}}{f'_{i,j,k}}=\frac{b'_{j,k'}}{b'_{j,k}}=\frac{f'_{i+1,j,k'}}{f'_{i+1,j,k}}$,
therefore by the same considerations: $\frac{b'_{j,k'}}{b'_{j,k}}$
does not depend on $j$. Let $a''_{i}=a'_{i}\cdot b_{1,1},\quad b''_{j}=\frac{b'_{j,1}}{b'_{1,1}}$
and $c''_{k}=\frac{b'_{1,k}}{b'_{1,1}}$, then for every $i,j,k\geq2$:
\begin{eqnarray*}
 &  & f'_{i,j,k}=a'_{i}b'_{j,k}=a'_{i}\cdot b_{1,1}(\frac{b'_{j,k}}{b'_{j,1}}\cdot\frac{b_{j,1}}{b_{1.1}})=a''_{i}\cdot b''_{j}\cdot c''_{k},
\end{eqnarray*}
so $F$ has the group form.
\end{proof}
\end{lem}
The following lemma gives us the definable field which generate the
field/twisted form. We use the notations of Definition \ref{def:frames}:
\begin{lem}
\label{lem:diagonal}Let $F=(f_{i,j,k})$ be a 3 indiscernible array
with $\mbox{\ensuremath{\alpha(m,n,p)=m+n+p-2}}$, and let $\mbox{\ensuremath{F'=(f_{i,j,j})_{i,j\in\mathbb{N}}}}$.
If $RM(\mathcal{C}_{\mathcal{I}_{3}})=0$ where $\mbox{\ensuremath{\mathcal{I}_{3}=\{\{f_{1,1,1},f_{2,1,1}\},\{f_{1,2,1},f_{2,2,1}\},\{f_{1,1,2},f_{2,1,2}\},\{f_{1,2,2}f_{2,2,2}\}\}}}$,
then $\alpha_{F'}(m,n)=m+2n-2$.
\begin{proof}
Let $A_{m,n}=\{f_{i,j,k}\,|\,1\leq i\leq m,\quad1\leq j,k\leq n\}$
be the $m\times n\times n$ cube. We  first show that $acl(A_{2,n})\subset acl(f_{i,j,j}\,|\,1\leq j\leq n,\quad1\leq i\leq2)$
$(\dagger)$ by induction. \\
\\
For $n=2$: if $acl(A_{2,2})\nsubseteq acl(f_{1,1,1},f_{1,2,2},f_{2,1,1},f_{2,2,2})$,
then $\mbox{\ensuremath{RM(f_{1,1,1},f_{1,2,2},f_{2,1,1},f_{2,2,2})<4}}$
so it must be equal to $3$. Thus by Lemma \ref{lem:Counter1}(2)
and Fact \ref{fct:2.2} there is a definable group $G$ and elements
$\mbox{\{\ensuremath{a_{i}}\},\{\ensuremath{b_{j}}\}\ensuremath{\subset}G}$
such that $acl(f_{i,j,j})=acl(a_{i}\cdot b_{j})$. In particular $g=\frac{a_{2}}{a_{1}}\in acl(f_{1,1,1},f_{2,1,1})\cap acl(f_{1,2,2}f_{2,2,2})$.\\
By indiscerniblity $g\in acl(f_{1,3,2},f_{2,3,2})\cap acl(f_{1,2,3},f_{2,2,2})\cap acl(f_{1,3,3},f_{2,3,3})$.
By using indiscerniblity one more time, we may retrieve an automorphism
$\sigma$ which sends the cube: $(f_{i,j,k}\,|\,1\leq i\leq2,\quad2\leq j,k\leq3)$
to the unit cube $(f_{i,j,k}\,|\,1\leq i,j,k\leq2)$, thus $\sigma(g)\in\bigcap_{l\in\mathcal{I}_{3}}acl(l)$,
in contradiction to the assumption that $Rm(\mathcal{C}_{\mathcal{I}_{3}})=0$.
\\
\\
If $n>2$, then by induction hypothesis: $acl(A_{2,n-1})\subset acl(f_{i,j,j}\,|\,1\leq j\leq n-1,\quad1\leq i\leq2)$.
Notice that: $\mbox{\ensuremath{tp(f_{i,j,k}\,|\,1\leq i,j,k\leq2)=tp(f_{i,j,k}\,|\,1\leq i\leq2,\quad j,k\in\{l,n\})}}$
( for every $1\leq l\leq n-1$ ), since the following holds as well
\begin{eqnarray*}
 &  & \{f_{i,j,k}\,|\,1\leq i\leq2\quad j,k\in\{l,n\}\}\subset acl(A_{2,n-1})\cup acl(f_{1,n,n},f_{2,n,n})\subset\\
 &  & acl(\{f_{i,j,j}\,|\,1\leq j\leq n\quad1\leq i\leq2\})
\end{eqnarray*}
We can use indiscerniblity to extend $(\dagger)$ to every $m\in\mathbb{N}$
and get:
\begin{eqnarray*}
 &  & acl(A_{m,n})\subset acl(f_{i,j,j}\,|\,1\leq j\leq n,\quad1\leq i\leq m).
\end{eqnarray*}
Thus $\mbox{\ensuremath{\alpha_{F'}}(m,n)=\ensuremath{\alpha_{F}}(m,n,n)=m+2n-2.}$
\end{proof}
\end{lem}
\begin{rem}
The same proof (just changing indices) will also work for the frames
$\mathcal{I}_{1}$ and $\mathcal{I}_{2}$. 

We are now ready to prove our Main Theorem.\end{rem}
\begin{thm}
Let $F=(f_{i,j,k})$ be a 3 indiscernible array and suppose that $\alpha_{F}(m,n,k)=m+n+k-2$,
then there are exactly three possibilities:
\begin{enumerate}
\item F has the group form.
\item F has the field form.
\item F has the twisted form.\end{enumerate}
\begin{proof}
By Lemma \ref{lem:dividing line} every form excludes the other ones,
so these are strictly distinct cases. In order to prove that $F$
must have one of the above forms we split to cases according to the
dividing line: \\
If for every frame of the unit cube $\mathcal{I}$ $Rm(\mathcal{C}_{\mathcal{I}})=1$,
then $F$ has the group form by Lemma \ref{lem:group form}. In the
other case, there is a frame $\mathcal{I}$ such that $Rm(\mathcal{C}_{\mathcal{I}})=0$
$(\dagger)$ .\\
 \\
Without loss of generality we may assume $\mathcal{I}=\mathcal{I}_{3}=\{\{f_{1,1,1},f_{2,1,1}\},\{f_{1,2,1},f_{2,2,1}\},\{f_{1,1,2},f_{2,1,2}\},\{f_{1,2,2}f_{2,2,2}\}\}$,
by Lemma \ref{lem:diagonal} and Fact \ref{fct:2.2} there is a definable
field $\mathbb{F}$ and elements $\mbox{\{\ensuremath{u_{i}}\},\{\ensuremath{v_{j}}\},\{\ensuremath{w_{k}}\}\ensuremath{\subset\mathbb{F}}}$
such that:\\
 $acl(f_{i,j,j})=acl(u_{i}v_{j}+w_{j})$, so we may use Lemma \ref{lem:Constructing lemma}
in order to get $\{u'_{i}\},\{v'_{j,k}\},\{w'_{j,k}\}\subset\mathbb{F}$
such that: $F'=(f'_{i,j,k})$ where $f'_{i,j,k}=u'_{i}v'_{j,k}+w'_{j,k}$
, is a 3 indiscernible array which is algebraically equivalent to
$F$, we may assume $v'_{1,1}=1$ and $w'_{1,1}=0$. Let $A_{1}=(f_{i,1,k}\,|\,1\leq i,k\leq2)$
and $\mbox{\ensuremath{A_{2}}=(\ensuremath{f_{i,1,k}}\,|\,3\ensuremath{\leq}i\ensuremath{\leq}4,\quad1\ensuremath{\leq}k\ensuremath{\leq}2)}$.
\\
Observe that $\mbox{\ensuremath{(w'_{1,2},v'_{1,2})\in acl(A_{1})\cap acl(A_{2})}}$,
hence:
\begin{eqnarray*}
 &  & RM(A_{1}A_{2})-RM(w'_{1,2},v'_{1,2})=RM(A_{1}A_{2}/(w'_{1,2},v'_{1,2}))\leq\\
 &  & RM(A_{1}/(w'_{1,2},v'_{1,2}))+RM(A_{2}/(w'_{1,2},v'_{1,2}))=2RM(A_{1})-2RM((w'_{1,2},v'_{1,2})),
\end{eqnarray*}
 and $RM(w'_{1,2},v'_{1,2})\leq6-5=1$. \\
On the other hand, $RM(A_{1})=3>2$, thus $RM(w'_{1,2},v'_{1,2})\geq1$,
so it must be equal to $1$. In the same way we get $RM(w'_{2,1},v'_{2,1})=1$.
We split into cases:
\begin{enumerate}
\item Assume that $w'_{1,2}\notin acl(v'_{1,2})$: then $RM(v'_{1,2})=RM(\frac{v'_{1,2}}{v'_{1,1}})=0$
and by indiscerniblity for every $j\in\mathbb{N}$: $RM(\frac{v'_{j,2}}{v'_{j,1}})=0$.
If also $RM(v'_{2,1})=0$, then $RM(v'_{2,2})=RM(\frac{v'_{2,2}}{v_{2,1}}\cdot v'_{2,1})=0$,
hence
\begin{eqnarray*}
 &  & acl(u_{2}-u_{1})=acl(v'_{2,2}(u_{2}-u_{1}))\in acl(f_{1,1,1},f_{2,1,1})\cap acl(f_{1,2,2}f_{2,2,2})
\end{eqnarray*}
 in contradiction to assumption $(\dagger)$. Thus: $RM(v'_{j,1})=1$
and $w'_{j,1}\in acl(v'_{j,1})$, and by Lemma \ref{lem:Field2} (1)+(2)
$F$ has the field form.
\item Assume that $w'_{1,2}\in acl(v'_{1,2})$ (in particular $RM(v'_{1,2})=1$):\\
If $w'_{2,1}\notin acl(v'_{2.1})$ then again by Lemma \ref{lem:Field2}(1)+(2)
$F$ has the field form.\\
If $w'_{2,1}\in acl(v'_{2,1})$ then by assumption $(\dagger)$ and
Lemma \ref{lem:Field2} (3), $F$ has the twisted form.\\

\end{enumerate}
\end{proof}
\end{thm}

\section{Application }

In this section we will see an application to Theorem \ref{thm:main}.
\\
\\
Let $k$ ba an algebraically closed field of characteristic zero and
$k\subset K$ be a universal domain of $k$. Every variety over $k$
will be identified with a definable subset in $K^{m}$ for some $m\in\mathbb{N}$
and each morphism between varieties will be identified with a definable
function between the corresponding definable sets (for more information
about the model theoretic framework see {[}5{]}, for the precise identification
see Remark 3.10 there).
\begin{defn}
Let $V_{1},V_{2},V_{3},U$ be irreducible curves over $k$ and let
$P:\,V_{1}\times V_{2}\times V_{3}\rightarrow U$ be a rational map.
\begin{enumerate}
\item $P$ is non-degenerate if for every generic independent elements $v_{i}\in V_{i}$,
the set: $\{v_{1},v_{2},v_{3},P(v_{1},v_{2},v_{3})\}$ is 3 independent.
\item $P$ has the group form, if there exists a one dimensional algebraic
group $(G,*)$, an affine variety $\overline{U}$ and rational maps:
$\mbox{\ensuremath{r_{i}:\,\ V_{i}\rightarrow G,\,\,\,\overline{P}:\,\ V_{1}\times V_{2}\times V_{3}\rightarrow\overline{U},\,\,\,\ q:\,\overline{U}\rightarrow G}}$,
and $\pi:\,\overline{U}\rightarrow U$ such that:
\begin{eqnarray*}
 &  & \pi\circ\overline{P}(x,y,z)=P(x,y,z)\mbox{ and }q\circ\overline{P}(x,y,z)=r_{1}(x)*r_{2}(y)*r_{3}(z).
\end{eqnarray*}

\item P has the field form if there exists $n\in\mathbb{N}$, an affine
variety $\overline{U}$ and rational maps:\\
 $\mbox{\ensuremath{\ensuremath{r_{i}:\ V_{i}\rightarrow K,\,\,\,\overline{P}:\ V_{1}\times V_{2}\times V_{3}\rightarrow\overline{U},\,\,\,q:\,\overline{U}\rightarrow K}}}$,
and $\pi:\,\overline{U}\rightarrow U$ such that:
\begin{eqnarray*}
 &  & \pi\circ\overline{P}(x,y,z)=P(x,y,z)\mbox{ and }q\circ\overline{P}(x_{1},x_{2},x_{3})=r_{i}(x_{i})\cdot(r_{j}(x_{j})+r_{l}(x_{l}))^{n},
\end{eqnarray*}
where $i\text{\ensuremath{\neq}}j\text{\ensuremath{\neq}}l\in\{1,2,3\}$.
\item P has the twisted form if there exists $n\in\mathbb{N}$, an affine
variety $\overline{U}$ and rational maps: \\
$\mbox{\ensuremath{r_{i}:\ V_{i}\rightarrow K,\,\,\overline{P}:\ V_{1}\times V_{2}\times V_{3}\rightarrow\overline{U},\,\,q:\,\overline{U}\rightarrow K}}$
, and $\pi:\,\overline{U}\rightarrow U$ such that:
\begin{eqnarray*}
 &  & \pi\circ\overline{P}(x,y,z)=P(x,y,z)\mbox{ and }q\circ\overline{P}(x,y,z)=\frac{r_{1}(x)+r_{2}(y)}{r_{2}(y)+r_{3}(z)}.
\end{eqnarray*}

\end{enumerate}
Our aim, is to prove the following theorem:\end{defn}
\begin{thm}
\label{thm:App}Let $V_{1},V_{2},V_{3},U$ be irreducible curves,
and let $P:\,V_{1}\times V_{2}\times V_{3}\rightarrow U$ be a non-degenerate
rational map such that:
\begin{eqnarray*}
 &  & T=\{(P(u_{i},v_{j},w_{k})\,|\,(i,j,k)\in\{0,1\}^{3})\,\,|\,u_{1},u_{2}\in V_{1},\,\,v_{1},v_{2}\in V_{2},\,\,w_{1},w_{2}\in V_{3}\}\subset U^{8}
\end{eqnarray*}
 has Zariski dimension 4, then P has either the group form, the field
form or the twisted form.\end{thm}
\begin{defn}
Let $V_{1},V_{2},V_{3},U$ be irreducible curves and let $P:\,V_{1}\times V_{2}\times V_{3}\rightarrow U$
be a rational map, we say that a 3 indiscernible array $F=(f_{i,j,k})$
is induced by $P$, if there exist an independent set: $\{u_{i},v_{j},w_{k}\,|\,u_{i}\in V_{1}\,\,v_{j}\in V_{2}\,\,w_{k}\in V_{3}\,\,i,j,k\in\mathbb{N}\}$
such that $f_{i,j,k}=P(u_{i},v_{j},w_{k})$.\end{defn}
\begin{lem}
\label{lem:bla1}Let $V_{1},V_{2},V_{3},U$ be irreducible curves,
$S:\,V_{1}\times V_{2}\times V_{3}\rightarrow U$ be a non-degenerate
rational map and $F=(f_{i,j,k})$ be a 3 indicernible array induced
by S. If 
\begin{eqnarray*}
 &  & T=\{(S(u_{i},v_{j},w_{k})\,|\,(i,j,k)\in\{0,1\}^{3})\,|\,u_{1},u_{2}\in V_{1},\,\,v_{1},v_{2}\in V_{2},\,\,w_{1},w_{2}\in V_{3}\}\subset U^{8}
\end{eqnarray*}
 has Zariski dimension 4, then $\alpha_{F}(m,n,p)=m+n+p-2$.
\begin{proof}
Let $\{u_{i},v_{j},w_{k}\,|\,u_{i}\in V_{1},\,\,v_{j}\in V_{2},\,\,w_{k}\in V_{3},\,\,i,j,k\in\mathbb{N}\}$
be an independent set such that $\mbox{\ensuremath{f_{i,j,k}}=S(\ensuremath{u_{i}},\ensuremath{v_{j}},\ensuremath{w_{k}})}$.
By our assumption 
\begin{eqnarray*}
 &  & RM(f_{1,j,k}\,|\,1\leq j,k\leq2)=RM(f_{j,1,k}\,|\,1\leq j,k\leq2)=RM(f_{j,k,1}\,|\,1\leq j,k\leq2)=3,
\end{eqnarray*}
 thus by Lemma \ref{lem:Counter1} $\alpha_{F}(m,n,p)=m+n+p-2$.\end{proof}
\begin{lem}
\label{lem:8.3}Let $V_{1},V_{2},V_{3},U$ be irreducible curves,
$P:\,V_{1}\times V_{2}\times V_{3}\rightarrow U$ a non-degenerate
rational map and $a_{i}\in V_{i}$ (where $1\leq i\leq3$) generic
independent elements. Suppose there exists a definable field $\mathbb{F}$,
and independent generic elements $u_{1},u_{2},u_{3}\in\mathbb{F},$
such that: $acl(a_{i})=acl(u_{i})$ for every $1\leq i\leq3$ and
$acl(P(a_{1},a_{2},a_{3}))=acl(\frac{u_{3}+u_{2}}{u_{2}+u_{1}})$
or $\mbox{\ensuremath{acl(P(a_{1},a_{2},a_{3}))=acl(u_{1}(u_{2}+u_{3}))}}$,
then: $u_{3}\in dcl(a_{3})$.
\begin{proof}
$ $

Case 1: If $acl(P(a_{1},a_{2},a_{3}))=acl(u_{1}(u_{2}+u_{3}))$, let
$u'_{3}$ be a conjugate of $u_{3}$ over $a_{3}$. Because $u_{3}\downarrow u_{2},u_{1}$,
then by uniqueness of the non forking extension: 
\begin{eqnarray*}
 &  & tp((u_{3},a_{3})/acl(u_{1},u_{2}))=tp((u'_{3},a_{3})/acl(u_{1},u_{2})).
\end{eqnarray*}
Hence there is an automorphism $\sigma$ which fixes $(a_{1},a_{2},a_{3},u_{2},u_{1})$
and sends $u_{3}$ to $u'_{3}$. Thus 
\begin{eqnarray*}
 &  & acl(u_{1}(u_{2}+u_{3}))=acl(u_{1}(u_{2}+u'_{3}))=acl(P(a_{1},a_{2},a_{3})),
\end{eqnarray*}
 and if we define $t=u_{1}(u_{2}+u_{3})$, then $u_{1}\cdot(u'_{3}-u_{3})\in acl(t)$.
However, $t\downarrow u_{1},u_{3}$ thus $\mbox{\ensuremath{u_{1}\cdot(u'_{3}-u_{3})\in acl(\emptyset)}}$,
but $\mbox{\ensuremath{(u'_{3}-u_{3})\in acl(u_{1})\cap acl(u_{3})=acl(\emptyset)}}$,
so we must have $u'_{3}=u_{3}$.\\
\\
 Case 2: If $acl(P(a_{1},a_{2},a_{3}))=acl(\frac{u_{1}+u_{2}}{u_{2}+u_{3}})$
then again, let $u'_{3}$ be a conjugate of $u_{1}$ over $a_{1}$.
Notice that: $acl(\frac{u_{3}+u_{2}}{u_{2}+u_{1}})=acl(\frac{u'_{3}+u_{2}}{u_{2}+u_{1}})$.
Let $t=\frac{u_{3}+u_{2}}{u_{2}+u_{1}}$ and $s=\frac{1}{u_{2}+u_{1}}$,
then $acl(t)=acl(t+s(u_{3}-u'_{3}))$, therefore, as in Case 1 we
must have: $u_{3}=u'_{3}$.

\end{proof}
\end{lem}
\end{lem}
We are working in an algebraically closed field of characteristic
zero, so if there are two irreducible varieties $V_{1},V_{2}$ and
generic elements $g\in V_{1},\,\,g'\in V_{2}$ such that $g\in dcl(g')$,
then there is a rational function $h:\,V_{1}\rightarrow V_{2}$ such
that $h(g)=g'$ (see {[}6{]} 4.9). We also use the fact that every
definable field $\mathbb{F}$ in K is definably isomorphic to K (see
{[}5{]} 4.13).
\begin{thm}
\label{thm:Rational-functions:}Let $V_{1},V_{2},V_{3},U$ be irreducible
curves, $P:\,V_{1}\times V_{2}\times V_{3}\rightarrow U$ a non-degenerate
rational map, and let $a_{i}\in V_{i}$ (where $1\leq i\leq3$) be
generic independent elements.
\begin{enumerate}
\item If there exist generic independent elements $\{u_{1},u_{2},u_{3}\}\subset K$
such that $acl(a_{i})=acl(u_{i})$ and $acl(P(a_{1},a_{2},a_{3}))=acl(\frac{u_{1}+u_{2}}{u_{2}+u_{3}})$,
then $P$ has the field form. 
\item If there exist generic independent elements $\{u_{1},u_{2},u_{3}\}\subset K$
such that $acl(a_{i})=acl(u_{i})$ and $acl(P(a_{1},a_{2},a_{3}))=acl(u_{1}(u_{2}+u_{3}))$,
then $P$ has the twisted form. \end{enumerate}
\begin{proof}
$ $
\begin{enumerate}
\item Observe that
\begin{eqnarray*}
 &  & acl(P(a_{1},a_{2},a_{3}))=acl(\frac{u_{3}+u_{2}}{u_{2}+u_{1}})=acl(\frac{u_{1}+u_{2}}{u_{2}+u_{3}})=acl(\frac{1}{\frac{u_{1}+u_{2}}{u_{2}+u_{3}}-1}=\frac{u_{2}+u_{3}}{u_{1}-u_{3}}),
\end{eqnarray*}
thus by Lemma \ref{lem:8.3} for every $1\leq i\leq3$: $acl(a_{i})=acl(u_{i})$,
so there are rational maps $r_{i}:\,V_{i}\rightarrow K$ such that
$r_{i}(a_{i})=u_{i}$. Choose $t\in k(P(a_{1},a_{2},a_{3}))^{alg}$
such that $k(a_{1},a_{2},a_{3})\cap k(t)^{acl}=k(t)$. Let $\overline{U}$
be an affine variety such that $k(t)=k(\overline{U})$. \\
On the one hand: $\frac{u_{3}+u_{2}}{u_{2}+u_{1}}\in dcl(t)$, on
the other hand: $P(a_{1},a_{2},a_{3})\in dcl(t)$ and $t\in dcl(a_{1},a_{2},a_{3})$.
Thus there exists rational maps: $\overline{P}:\,V_{1}\times V_{2}\times V_{3}\rightarrow\overline{U}$,
$q:\,\overline{U}\rightarrow K$ and $\pi:\,\overline{U}\rightarrow U$
such that:\\
$\pi\circ\overline{P}(x,y,z)=P(x,y,z)$ and $q\circ\overline{P}(x,y,z)=\frac{r_{1}(x)+r_{2}(y)}{r_{2}(y)+r_{3}(z)}$.\\

\item By lemma \ref{lem:8.3} there are rational maps $r_{i}:\,V_{i}\rightarrow K$
such that $r_{i}(a_{i})=u_{i}$ where $2\leq i\leq3$.\\
Let $u'_{1}$ be a conjugate of $u_{1}$ over $a_{1}$, then, as in
the previous lemma, there is an automorphism $\sigma$ which fixes
$a_{1},a_{2},a_{3},u_{2},u_{3}$ and sends $u_{1}$ to $u'_{1}$.
So $acl(u_{1}(u_{2}+u_{3})=acl(u'_{1}(u_{2}+u_{3}))$. If we take
$t=u_{1}(u_{2}+u_{3})$ then $\sigma(t)=t\cdot\frac{u'_{1}}{u_{1}}$
and $acl(t\cdot\frac{u'_{1}}{u_{1}})=acl(t)$. Hence, $\frac{u'_{1}}{u_{1}}\in acl(t)\cap acl(u_{1})=acl(\emptyset)$
and $\sigma^{i}(t)=t\cdot(\frac{u'_{1}}{u_{1}})^{n}\in acl(a_{1},a_{2},a_{3})$
, so there must be $n\in\mathbb{N}$ such that $(\frac{u'_{1}}{u_{1}})^{n}=1$.\\
Fix $N\in\mathbb{N}$ such that $(\frac{u'}{u''})^{N}=1$ for every
two conjugates $u',u''$ of $u_{1}$ over $a_{1}$. $u_{1}^{N}\in dcl(a_{1})$
so there is a rational map $r_{1}:\,V_{1}\rightarrow\mathbb{F}$ such
that $r_{1}(a_{1})=u_{1}^{N}$. In conclusion: 
\begin{eqnarray*}
 &  & acl(P(a_{1},a_{2},a_{3}))=acl(r_{1}(a_{1})\cdot(r_{2}(a_{2})+r_{3}(a_{3}))^{N})
\end{eqnarray*}
\newpage{}and as in Part 1 of this theorem we can find an affine
variety $\overline{U}$, and rational maps $\overline{P},\pi,q$ such
that: $\pi\circ\overline{P}(x,y,z)=P(x,y,z)$ and $q\circ\overline{P}(x,y,z)=r_{1}(x)\cdot(r_{2}(y)+r_{3}(z))^{N}$.
\end{enumerate}
\end{proof}
\begin{thm}
\label{thm:group}Let $V_{1},V_{2},V_{3},U$ be irreducible curves,
$P:\,V_{1}\times V_{2}\times V_{3}\rightarrow U$ a non-degenerate
rational map, $a_{i}\in V_{i}$ (where $1\leq i\leq3$) be generic
independent elements, and $(G,*)$ be a one dimensional algebraic
group. If there exist generic independent elements $\{u_{1},u_{2},u_{3}\}\subset G$
such that $acl(a_{i})=acl(u_{i})$ and $acl(P(a_{1},a_{2},a_{3}))=acl(u_{1}*u_{2}*u_{3})$
then $P$ has the group form.
\begin{proof}
As in the proof of Theorem \ref{thm:Rational-functions:} we can find
$N\in\mathbb{N}$ such that $u_{i}^{N}\in dcl(a_{i})$  Thus there
are rational maps $r_{i}:\,V_{i}\rightarrow G$ such that $r_{i}(a_{i})=u_{i}$
where $1\leq i\leq3$. Furthermore, we can find an affine variety
$\overline{U}$ and rational maps $\overline{P},\pi,q$ such that:
$\mbox{\ensuremath{\pi\circ\overline{P}(x,y,z)=P(x,y,z)}}$ and $\mbox{\ensuremath{q\circ\overline{P}(x,y,z)=r_{1}(x)*r_{2}(y)*r_{3}(z)}}$.
\end{proof}
\end{thm}
\end{thm}
We are now finally ready to prove Theorem \ref{thm:App}:
\begin{proof}
$ $ \\
Let $\{u_{i},v_{j},w_{k}\,|\,u_{i}\in V_{1}\,\,v_{j}\in V_{2}\,\,w_{k}\in V_{3}\,\,i,j,k\in\mathbb{N}\}$
be an independent set such that: \\
$f_{i,j,k}=P(u_{i},v_{j},w_{k})$, and set $F=(f_{i,j,k})$. By Lemma
\ref{lem:bla1} $\alpha_{F}(m,n,p)=m+n+p-2$, so by Theorem \ref{thm:main}
there exists a definable algebraic structure $\mathbb{G}$ (where
$\mathbb{G}$ is a one dimensional group in the case of group form
and field in the other cases), and independent elements: $\mbox{\ensuremath{\{a_{i},b_{j},c_{k}\}\subset\mathbb{G}}}$
such that: $acl(P(u_{i},v_{j},w_{k}))=acl(h(a_{i},b_{j},c_{k}))$,
where 
\begin{eqnarray*}
 &  & h(x,y,z)=x*y*z\mbox{ or }h(x,y,z)=x\cdot(y+z)\mbox{ or }h(x,y,z)=\frac{x+y}{y+z}
\end{eqnarray*}
 (the operations in $h$ are the operations of $\mathbb{G}$). In
either cases, if $k'=k(a_{1},b_{1},c_{1},u_{1},v_{1},w_{1})^{acl}$
then:\\
$acl(u_{i}/k')=acl(a_{i}/k')=acl(f_{i,1,1}/k'),\,\,acl(v_{j}/k')=acl(b_{j}/k')=acl(f_{1,j,1}/k')$
and\\
 $acl(w_{k}/k')=acl(c_{k}/k')=acl(f_{1,1,k}/k')$. By Theorems \ref{thm:Rational-functions:}
and \ref{thm:group} there exist an affine variety $\overline{U}$
and rational maps $\overline{P},\pi,q,r_{1},r_{2},r_{3}$ over $k'$
such that: 
\begin{eqnarray*}
 &  & \pi\circ\overline{P}(x,y,z)=P(x,y,z)\mbox{ and }q\circ\overline{P}(x,y,z)=h'(r_{1}(x),r_{2}(y),r_{3}(z)),
\end{eqnarray*}
(where $h'(x,y,z)=x*y*z$ in the group form, $\mbox{\ensuremath{h'(x,y,z)=x\cdot(y+z)^{n}}}$
in the field form and $\mbox{\ensuremath{h'(x,y,z)=\frac{x+y}{y+z}}}$
in the twisted form). However, $k$ is an elementary submodel of $k'$
, thus $\overline{U},\overline{P},\pi,q,r_{1},r_{2},r_{3}$ are definable
over $k$ as desired.
\end{proof}
One may use Theorem \ref{thm:App} to ``decompose'' (in a manner
we will details shortly) families of rational functions.
\begin{defn}
\label{Def:2-dec}Let $P(x,y,z)$ be a rational function over an algebraically
closed field $k$ of characteristic zero. We say that $P$ is 2-decomposed
if there exist rational functions $h_{1},h_{2},h_{3},v_{1},v_{2},v_{3}\in k(x,y)$
and $u_{1},u_{2},u_{3}\in k(x)$ such that:
\begin{enumerate}
\item $P(x,y,z)=u_{1}(h_{1}(y,x)\nabla^{1}v_{1}(z,x))$
\item $P(x,y,z)=u_{2}(h_{2}(x,y)\nabla^{2}v_{2}(z,y))$
\item $P(x,y,z)=u_{3}(h_{3}(x,z)\nabla^{3}v_{3}(y,z))$
\end{enumerate}

Where $\nabla^{i}\in\{\cdot,+\}$ for every $1\leq i\leq3$.

\end{defn}
\begin{cor}
\label{cor:Let--b}Let $P(x,y,z)$ be a rational function over an
algebraically closed field $k$ of characteristic zero. If $P$ is
a 2-decomposed rational function, then either there exist $n\in\mathbb{N}$
and rational functions $r_{1},r_{2},r_{3},q\in k(x)$ and $\overline{P}\in k(x,y,z)$
such that one of the following holds:
\begin{enumerate}
\item $q\circ\overline{P}(x_{1},x_{2},x_{3})=r_{i}(x_{i})(r_{j}(x_{j})+r(x_{l}))^{n}$
(for some $i\neq j\neq l\in\{1,2,3\})$ \\
and $\pi\circ\overline{P}(x,y,z)=P(x,y,z)$. 
\item $q\circ\overline{P}(x,y,z)=\frac{r_{1}(x)+r_{2}(y)}{r_{2}(y)+r_{3}(z)}$
and $\pi\circ\overline{P}(x,y,z)=P(x,y,z)$. 
\end{enumerate}

or there exist a one dimension algebraic group $(G,*)$, rational
maps \textup{$r_{1},r_{2},r_{3},q$ and $\overline{P}$ such that:}\\
\textup{$q\circ\overline{P}(x,y,z)=r_{1}(x)*r_{2}(y)*r_{3}(z)$.}
\begin{proof}
If $P$ is 2-decomposed then: 
\begin{eqnarray*}
 &  & T=\{(S(u_{i},v_{j},w_{k})\,|\,(i,j,k)\in\{0,1\}^{3})\,|\,u_{1},u_{2}\in V_{1}\,\,v_{1},v_{2}\in V_{2}\,\,w_{1},w_{2}\in V_{3}\}\subset U^{8}
\end{eqnarray*}
 has Zariski dimension 4, therefore, we can use Theorem \ref{thm:App}.
\end{proof}
\end{cor}
We end with an open question.
\begin{problem}
If $P(x,y,z)$ is a 2-decomposed polynomial, what more can be said?
\end{problem}
We expect a decomposition in the spirit of {[}1{]}:\\
There exist $n\in\mathbb{N}$ and polynomials $r_{1},r_{2},r_{3},u\in k[x]$,
such that $\mbox{\ensuremath{P(x_{1},x_{2},x_{3})=u(r_{i}(x_{i})\cdot(r_{j}(x_{j})+r_{l}(x_{l}))^{n})}}$,
where $i\neq j\neq l\in\{1,2,3\}$ or $\mbox{\ensuremath{P(x_{1},x_{2},x_{3})=u(r_{1}(x_{1})\nabla^{4}r_{2}(x_{2})\nabla^{4}r_{3}(x_{3}))}}$,
$\nabla^{4}\in\{\cdot,+\}$.

\section*{References}

{[}1{]} Tao, Terence. Expanding Polynomials Over Finite Fields Of
Large Characteristic, And A Regu-\\
$\mbox{}\quad\,\,\,$larity Lemma For Definable Sets arXiv:1211.2894\\
{[}2{]} Hrushovski, Ehud. Expanding comments \\
{[}3{]} Pillay, Anand. Geometric Stability Theory. Oxford: Clarendon
Press, 1996. Print.\\
{[}4{]} Hrushovski, Ehud, and Boris Zilber. Zariski Geometries. Bulletin
of the American Mathematical \\
$\mbox{}\quad\,\,\,$Society 28.2 (1993): 315-324. Web.\\
{[}5{]} Bouscaren, Elisabeth. Model Theory And Algebraic Geometry.
Berlin: Springer, 1998. Print.\\
{[}6{]} Poizat, Bruno. Stable Groups. Providence, R.I.: American Mathematical
Society, 2001. Print.
\end{document}